\documentclass[10pt,a4paper,leqno]{article}
\usepackage{amsmath,amssymb,graphicx}
\usepackage{latexsym}
\usepackage[english,francais]{babel}
\usepackage[latin1]{inputenc}
\usepackage{url}
\def\boxit#1#2{\vbox{\hrule\hbox{\vrule\vbox spread#1{\vfil\hbox spread#1{\hfil#2\hfil}\vfil}
\vrule}\hrule}}

\long\def\emboite#1#2{\setbox4=\vbox{\hsize #1\noindent \strut #2}
  \boxit{8pt}{\box4}}

\def\R{\mathbb R}

\def\C{\mathbb C}
\def\N{\mathbb N}

\newcommand{\bpmat}{\begin{pmatrix}}
\newcommand{\epmat}{\end{pmatrix}}

\ifx\optionkeymacros\undefined\else\endinput\fi

\catcode`\=\active\def{\aa}         % option a
\catcode`\º=\active\defº{\int}        % option b (math mode) 
\catcode`\=\active\def{\c c}        % option c
\catcode`\¶=\active\def¶{\partial}    % option d (math mode)
\catcode`\Ä=\active\defÄ{\oint}       % option f (math mode) ?
\catcode`\Æ=\active\defÆ{\triangle}   % option j (math mode)
\catcode`\Â=\active\defÂ{\neg}        % option l (math mode)
\catcode`\µ=\active\defµ{\mu}         % option m (math mode)
\catcode`\¿=\active\def¿{\o}          % option o
\catcode`\¹=\active\def¹{\pi}         % option p (math mode w/ arg.)
\catcode`\Ï=\active\defÏ{\oe}         % option q 
\catcode`\§=\active\def§{\ss}         % option s 
\catcode`\ =\active\def {\dagger}     % option t  (math mode)
\catcode`\Ã=\active\defÃ{\sqrt}       % option v (math mode w/ arg.)
\catcode`\·=\active\def·{\Sigma}      % option w (math mode)
\catcode`\Å=\active\defÅ{\approx}     % option x (math mode)
\catcode`\½=\active\def½{\Omega}      % option z (math mode)
\catcode`\£=\active\def£{{\it\$}}     % option 3 ($ from italic font)
\catcode`\°=\active\def°{\infty}      % option 5 (math mode)
\catcode`\¤=\active\def¤{\S}          % option 6 
\catcode`\¦=\active\def¦{\P}          % option 7
\catcode`\¥=\active\def¥{\bullet}     % option 8 
\catcode`\»=\active\def»{\b a}        % option 9
\catcode`\¼=\active\def¼{\b o}        % option 0
\catcode`\­=\active\def­{\not=}       % option = (math mode)
\catcode`\²=\active\def²{\leq}        % option , (math mode)
\catcode`\³=\active\def³{\geq}        % option . (math mode)
\catcode`\Ö=\active\defÖ{\div}        % option / (math mode)
\catcode`\É=\active\defÉ{\dots}       % option ; 
\catcode`\¾=\active\def¾{\ae}         % option '
\catcode`\Ç=\active\defÇ{\ll}         % option \ (math mode)
\catcode`\Ò=\active\defÒ{``}          % option [
\catcode`\=\active\def{\AA}         % shift-option A
\catcode`\=\active\def{\c C}        % shift-option C
\catcode`\¯=\active\def¯{\O}          % shift-option O
\catcode`\¸=\active\def¸{\Pi}         % shift-option P (math mode)
\catcode`\Î=\active\defÎ{\OE}         % shift-option Q
\catcode`\®=\active\def®{\AE}         % shift-option '
\catcode`\×=\active\def×{\diamond}    % shift-option V (math mode)
\catcode`\¡=\active\def¡{\accent'27}  % shift-option 8
\catcode`\Ó=\active\defÓ{''}          % shift-option [
\catcode`\±=\active\def±{\pm}         % shift-option = (math mode)
\catcode`\È=\active\defÈ{\gg}         % shift-option \ (math mode)
\catcode`\À=\active\defÀ{?`}          % shift-option / 
\catcode`\Ð=\active\defÐ{--}          % option - (en-dash)
\catcode`\Ñ=\active\defÑ{---}         % shift-option - (em-dash)

\catcode`\=\active\def{\"a}        % option u, then  a
\catcode`\=\active\def{\"e}        % option u, then  e
\catcode`\=\active\def{\"\i}       % option u, then  i
\catcode`\=\active\def{\"o}        % option u, then  o
\catcode`\=\active\def{\"u}        % option u, then  u
\catcode`\Ø=\active\defØ{\"y}        % option u, then  y
\catcode`\=\active\def{\"A}        % option u, then  A
\catcode`\=\active\def{\"O}        % option u, then  O
\catcode`\=\active\def{\"U}        % option u, then  U
\catcode`\=\active\def{\'a}        % option e, then  a
\catcode`\=\active\def{\'e}        % option e, then  e
\catcode`\=\active\def{\'\i}       % option e, then  i
\catcode`\=\active\def{\'o}        % option e, then  o
\catcode`\=\active\def{\'u}        % option e, then  u
\catcode`\=\active\def{\'E}        % option e, then  E
\catcode`\=\active\def{\`a}        % option `, then  a
\catcode`\=\active\def{\`e}        % option `, then  e
\catcode`\=\active\def{\`\i}       % option `, then  i
\catcode`\=\active\def{\`o}        % option `, then  o
\catcode`\=\active\def{\`u}        % option `, then  u
\catcode`\Ë=\active\defË{\`A}        % option `, then  A
\catcode`\=\active\def{\~a}        % option n, then  a
\catcode`\=\active\def{\~n}        % option n, then  n
\catcode`\=\active\def{\~o}        % option n, then  o
\catcode`\Ì=\active\defÌ{\~A}        % option n, then  A
\catcode`\=\active\def{\~N}        % option n, then  N
\catcode`\Í=\active\defÍ{\~O}        % option n, then  O
\catcode`\=\active\def{\^a}        % option i, then  a
\catcode`\=\active\def{\^e}        % option i, then  e
\catcode`\=\active\def{\^\i}       % option i, then  i
\catcode`\=\active\def{\^o}        % option i, then  o
\catcode`\=\active\def{\^u}        % option i, then  u

\let\optionkeymacros\null

\font\std=cmr10 at 10pt

\font\biengras=cmbx10 at 15pt

\font\ppbf=cmbx10 at 8pt

\def\beauno#1={\oldstyle at 14pt #1}

\def\cqfd{\hfill $\square$
\vskip0,1cm\kern0,3cm\noindent}

\newcount\num 
\newcount\sousnum 
\def\nu{\noindent
\sousnum=0\advance\num by 1
\if\the\num1{\bf\the\num.\ }
\else{\ \bf\par\noindent\the\num.\ }
\fi
}

\def\snu{\par\noindent
\advance\sousnum by 1
{\bf\the\num.\the\sousnum.\ }}

\newcount\temoin 
\def\exo{\num=0\advance\temoin by 1
\quad

\medskip\noindent{\bf E{\ppbf XERCICE}\
\the\temoin\quad\std\hrulefill\par\noindent}}

\def\hv{\par\kern-0.15cm\noindent
\hrulefill\ \par\bigskip\noindent}

\def\nspace#1{\vskip0.1cm\kern-#1cm}

\begin{document}
\def\labelitemi{$\bullet$}
\thispagestyle{empty}

%% en-tete 
%\fancyhead[RO]{\small\slshape\leftmark\hfill P -- \thepage}%%TITRE COURT%%

\emboite{\hsize}{
  \vskip7pt
  \centerline {\textsf{\biengras{Transversalit quantitative en gomtrie symplectique :}}}
  \vskip7pt
  \centerline{\biengras{sous-varits et hypersurfaces}}
\vskip7pt
 \medskip
 
 }
% \quad \hfill{\std{Octobre} 2010}}
 % \\\\
  {\it 
    \noindent
 % \hrulefill
  }

\bigskip

\section{Enonc du thorme}
Un thorme de Donaldson et Auroux \cite{Do96} \cite{Au97} permet de perturber une section approximativement holomorphe d'un fibr trs positif afin de rendre son annulation quantitativement transversale. Ce thorme admet une variante relative  une sous-varit, c'est--dire qu'on peut rendre quantitativement transversale l'annulation de la restriction de la section  une sous-varit donne.

Afin d'noncer prcisment ce rsultat, donnons plusieurs dfinitions.
Soit une application linaire $u$ entre deux espaces euclidiens $E$ et $F$.
On appelle module d'injectivit de $u$ le nombre :
\begin{eqnarray*}
\mbox{MI}(u) & = & \mbox{min} \left\| u(v) \right\| 
\end{eqnarray*}
o $v$ dcrit la sphre unit de l'espace de dpart $E$. De faon duale,
on appelle module de surjectivit de $u$ le nombre :
\begin{eqnarray*}
\mbox{MS}(u) & = & \mbox{min} \left\| \lambda \circ u \right\| 
\end{eqnarray*}
o $\lambda$ dcrit la sphre unit de l'espace dual $F^*$ de l'espace d'arrive (autrement dit :
$ \mbox{MS}(u) = \mbox{MI}(u^*) $, o $u^*$ dsigne l'adjoint de $u$).

Soient une application lisse $f$ d'un ouvert $U \subset E$ vers $F$ et $x$ un point de $U$. On appelle module de transversalit de $f$ en $x$ le nombre :
\begin{eqnarray*}
\mbox{MT}(f, \, x) & = & \mbox{max} \left\{ \left\| f(x) \right\| ,\; \mbox{MS}(df(x))\right\}.
\end{eqnarray*}

Etendons la dfinition du module de transversalit aux sections :
soit un fibr vectoriel $V$ sur une varit $X$, muni d'une connexion $\nabla$. On suppose la base et la fibre munies de mtriques.
Alors on appelle module de transversalit d'une section $s$ de $V$ en un point $x$ de $X$ le nombre : 
\begin{eqnarray*}
\mbox{MT}(s, \, x) & = & \mbox{max} \left\{ \left\| s(x) \right\| ,\; \mbox{MS}(\nabla s(x))\right\}
\end{eqnarray*}
et on appelle module (global) de transversalit de $s$ le nombre :
$\mbox{MT}(s)  =  \mbox{inf} \; \mbox{MT}(s,\, x).$

Si $Y$ dsigne une sous-varit de $X$, on appellera module de transversalit de $s$ le long de $Y$ le module de transversalit de la restriction de $s$  $Y$.
Enfin, il existe de cette notion de module de transversalit une version pondre : le module de transversalit de poids $(a,\, b)$, avec $a$ et $b>0$, est dfini par :
\begin{eqnarray*}
\mbox{MT}(s, \, x;\; a,\, b) & = & \mbox{max} \left\{ a\times \left\| s(x) \right\| ,\; b \times \mbox{MS}(\nabla s(x))\right\}
\end{eqnarray*}
et les modules pondrs global et relatif  une sous-varit se dfinissent de faon analogue.
\\
\\
\indent
Rappelons aussi quelques dfinitions de gomtrie symplectique (ou de gomtrie presque khlrienne).
Une varit presque khlrienne est une varit de dimension (relle) paire munie d'une forme symplectique $\omega$, d'une structure presque complexe $J$ et d'une mtrique riemannienne $g$, ces trois donnes vrifiant :
\begin{eqnarray*}
g (v,\, w) & = & \omega(v,\, Jw).
\end{eqnarray*}
pour tous champs de vecteurs $v$ et $w$.

Une pr-quantification d'une varit presque khlrienne $(X,\, \omega,\,  J,\,  g)$ est un fibr en droites hermitiennes de base $X$, muni d'une connexion unitaire de courbure $F$ relie  la forme symplectique par la formule :
\begin{eqnarray*}
F & = & -i2\pi \omega.
\end{eqnarray*}
\\
\\
\indent
Le contexte de la thorie de Donaldson-Auroux est une varit presque khlrienne $X$ munie d'une pr-quantification $L$. Si $r$ et $k$ dsignent deux entiers $>0$, on notera $L^k$ la $k-$me puissance tensorielle de $L$ (c'est un fibr en droites) et on notera $rL^k$ son multiple de rang $r$. Certaines estimes, omniprsentes dans cette thorie, donnent lieu aux dfinitions suivantes.

Soient un rel $K>0$ et trois entiers $m_{max} \geq 0$, $r>0$ et $k>0$.
Une section $s$ du fibr hermitien $rL^k$ sera dite $(K,m_{max})-$contrle si elle vrifie :
\begin{eqnarray*}
\left\| \nabla^m s \right\| & \leq & K\times k^{\frac{m}{2}}
\end{eqnarray*}
pour tout entier $m$ compris entre $0$ et $m_{max} $.
Elle sera dite $(K,m_{max} )-$approximativement holomorphe si elle vrifie :
\begin{eqnarray*}
\left\| \nabla^m \overline{\partial} s \right\| & \leq & K \times k^{\frac{m}{2}}
\end{eqnarray*}
pour tout entier $m$ compris entre $0$ et $m_{max} $.

Prcisons que $\overline{\partial}$ dsigne la partie anti-linaire de $\nabla$ et que toutes les mtriques et connexions qui interviennent dans 
les deux dfinitions prcdentes sont dduites des mtriques de $X$ et de $L$, de la connexion de $L$ et de la connexion de Levi-Civita. 
\\
\\
\indent
Nous pouvons maintenant noncer le principal rsultat du prsent article :
\\
\\
\indent
{\bf Thorme de Donaldson-Auroux relatif ( une sous-varit).}
Soient une varit presque khlrienne $X$, une pr-quantification $L$, une sous-varit compacte $Y\subset X$, deux rels $\varepsilon >0$ et $K>0$ et deux entiers $m_{max}  \geq 0$ et $r>0$. Alors il existe existe un rel $\eta>0$ et un entier $k_0>0$ tels que, pour tout entier $k \geq k_0$ et toute section $(K,2)-$contrle et $(K,1)-$approximativement holomorphe $s$ de $rL^k$, il existe une section $(\varepsilon,m_{max} )-$contrle et $(\varepsilon,m_{max} )-$approximativement holomorphe $t$ de $rL^k$ telle que le module de transversalit de poids $(1,\, k^{-\frac{1}{2}})$ de $s+t$ le long de $Y$ soit suprieur  $\eta$.
\\
\\
\indent
Quoiqu'assez longue, la dmonstration de ce thorme n'a rien de mystrieux : il s'agit d'une adaptation de la preuve originale de Donaldson. 
Le principal intrt de la version relative du thorme de Donaldson-Auroux est de retrouver comme corollaires les deux rsultats suivants :

-Le thorme de D. Auroux sur l'unicit  isotopie prs des sections quantitativement transversales de Donaldson \cite{Au97}.
La dmonstration originale reposait sur un argument de \og grande composante connexe \fg.
Celle propose ici en diffre donc sur ce point.

-Le thorme de A. Ibort, D. Mart\'{i}nez-Torres et F. Presas \cite{IbMaPr00}, analogue en gomtrie de contact du thorme de Donaldson-Auroux.
\\
\\
\indent
En fait, ces deux rsultats sont des consquences du thorme de Donaldson-Auroux relatif aux {\it hypersurfaces relles}.
Dans le cas des hypersurfaces relles, le point important est qu'une section (approximativement) holomorphe 
qui s'annule transversalement (avec une estime) s'annulera transversalement {\it le long de la distribution de Levi}.
\\
\\
\indent
{\bf Plan de l'article.}

La premire partie de l'article (chapitres \ref{chapimp} et \ref{chapau}) est consacre aux corollaires du thorme de Donaldson-Auroux relatif
et la seconde partie (chapitres \ref{chapvv} et suivants) est consacre  la dmonstration de ce thorme.

Le chapitre \ref{chapimp} traite le cas de la transversalisation le long d'une hypersurface relle puis
en dduit le thorme d'Ibort, Mart\'{i}nez-Torres et Presas.
De mme, dans le chapitre \ref{chapau}, le thorme d'unicit d'Auroux  est dduit du rsultat sur les hypersurfaces et ce chapitre
contient aussi un contre-exemple simple qui prouve qu'il n'y a pas unicit dans le cas relatif. 

Le chapitre \ref{chapvv} expose la thorie des variations de Vitushkin et l'applique  l'tude de la complexit
des hypersurfaces algbriques (d'aprs \cite{CoYo04}).  
L'objet du chapitre \ref{chapbert} est un thorme de Bertini quantitatif (inspir de \cite{Do99}).
Par souci de compltude, on donne, dans les chapitres suivants, la fin de la dmonstration du thorme de Donaldson-Auroux relatif
qui est trs semblable  la dmonstration classique de \cite{Do96}, \cite{Au97} et \cite{Do99}.
\\
\\
\indent
L'auteur remercie Emmanuel Giroux et Bruno Svennec pour leur soutien et pour de nombreuses discussions.
Merci aussi  Nicolas Dutertre et  Stphane Rigat pour des indications bibliographiques.

\section{Application  la gomtrie de contact}\label{chapimp}

Ce chapitre permet, en ce qui concerne les techniques de transversalit quantitative, de regarder la thorie de contact comme une application
de la thorie symplectique relative.

\subsection{Hyperplan rel et hyperplan complexe}\label{h et h}
Lemme. 
Soient deux espaces hermitiens $E$ et $F$ de dimensions finies, une application $u:E\mapsto F$ et un hyperplan rel $H \subset E$.
Notons $K$ l'hyperplan complexe inclus dans $H$ (autrement dit $K=H\cap iH$).
On dsignera par $u_{H}$ et $u_{K}$ les restrictions respectives
de $u$  $H$ et  $K$. Alors :

(i) si $u$ est  $\C -$linaire, elle vrifiera l'identit :
$$
\mbox{MS}(u_{K}) = \mbox{MS}(u_{H}),
$$
\indent
(ii) si $u$ est $\R -$linaire, elle vrifiera l'encadrement :
$$
\mbox{MS}(u_{H}) - 2 \|u^{0,1} \|\leq \mbox{MS}(u_{K}) \leq \mbox{MS}(u_{H})
$$
o $u^{0,1}$ dsigne la partie $\C -$antilinaire de $u$.
\\
\\
\indent
Dmonstration.
Dans les cas (i) et (ii), l'ingalit $\mbox{MS}(u_{K}) \leq \mbox{MS}(u_{H})$ est vrifie car $K$ est inclus dans $H$.

Cas (i). 
Ecrivons la dcomposition : 
$$
H = \R v_0 \oplus K
$$
o $v_0$ dsigne un vecteur orthogonal  $K$.
Par dfinition du module de surjectivit, il existe une forme linaire relle $\lambda: F \rightarrow \R$
de norme $1$ vrifiant :
$$
\| \lambda \circ u_{K} \| =\mbox{MS}(u_{K}).
$$ 
Le thorme des valeurs intermdiaires fournit un rel $\theta$ qui vrifie :
$$
\lambda ( e^{i\theta} u(v_0) )  = 0.
$$
Notons $\lambda_{\theta}$ la forme linaire relle dfinie sur $F$ par :
$$
\lambda_{\theta} (v) = \lambda (e^{i\theta} v).
$$
Sa norme est gale  celle de $\lambda$, c'est--dire vaut $1$.
Par dfinition du module de surjectivit, on en dduit l'ingalit :
$$
\mbox{MS}(u_{H})  \leq  \| \lambda_{\theta} \circ u_{H}\| .
$$
Comme $\lambda_{\theta} \circ u$ est nulle sur $\R v_0$, le membre de droite est gal  $\| \lambda_{\theta} \circ u_{K}\|$.
En appliquant la $\C-$linarit de $u_{K}$ puis la dfinition de $\lambda$, on peut crire :
$$
\| \lambda_{\theta} \circ u_{K}\| = \| \lambda \circ u_{K}\| = \mbox{MS}(u_{K}).
$$
On obtient donc l'ingalit suivante :
$$ 
\mbox{MS}(u_{H}) \leq \mbox{MS}(u_{K}).$$
Le cas (i) du lemme est dmontr par double ingalit.
\\
\\
\indent
Cas (ii). L'application $u$ est dsormais $\R -$linaire. 
Alors :
\begin{eqnarray*}
\mbox{MS}(u_{K}) & \geq &  \mbox{MS}((u^{1,0})_K) - \| (u^{0,1})_K \|
\\ & & \mbox{(car le module de surjectivit est $1 -$lipschitzien)}
\\ & = & \mbox{MS}((u^{1,0})_H) - \| (u^{0,1})_K \|
\\ && \mbox{(d'aprs (i) appliqu  $u^{1,0}$)}
\\ & \geq &  \mbox{MS}(u_{H}) - \| (u^{0,1})_H \| - \| (u^{0,1})_K \|
\\ && \mbox{(car le module de surjectivit est $1 -$lipschitzien)}
\\ & \geq &  \mbox{MS}(u_{H}) - 2 \| u^{0,1} \|
.\end{eqnarray*}
Le lemme est dmontr.

\subsection{Cas d'une hypersurface relle : transversalit dans la direction de Levi}\label{Levi}

Thorme de transversalisation le long d'une hypersurface relle.

Sous les hypothses du thorme de Donalsdson-Auroux relatif, on suppose de plus que $Y$ est une hypersurface relle.
Alors il existe existe (comme d'habitude) un rel $\eta>0$ et un entier $k_0>0$ tels que, 
pour tout entier $k \geq k_0$ et toute section $(K,2)-$contrle et $(K,1)-$approximativement holomorphe $s$ de $rL^k$, 
il existe une section $t$ de $rL^k$ qui est comme d'habitude $(\varepsilon,m_{max} )-$contrle et $(\varepsilon,m_{max} )-$approximativement holomorphe 
mais qui vrifie, pour tout $x\in Y$, l'estimation suivante :
$$
\mbox{max} \left\{\left\| (s+t)(x) \right\| ,\; k^{-\frac{1}{2}} \times \mbox{MS}(\nabla (s+t)(x)_{\mbox{Levi}})\right\}  \geq \eta
$$
o $\nabla (s+t)(x)_{\mbox{Levi}}$ dsigne la restriction de $\nabla (s+t)(x)$  l'hyperplan (complexe) de Levi $T_xY \cap J_x(T_xY)$.
\\
\\
\indent
Dmonstration.
Le thorme de Donaldson-Auroux relatif fournit une section $t$ satisfaisant les majorations demandes
et la minoration suivante :
$$
\mbox{max} \left\{\left\| (s+t)(x) \right\| ,\; k^{-\frac{1}{2}} \times \mbox{MS}(\nabla (s+t)(x)_{T_x Y})\right\}  \geq \eta_1
$$
pour un certain rel $\eta_1 > 0$. Posons $\eta = \frac{\eta_1}{2}$.
Dans le cas particulier d'une hypersurface relle, le lemme (\ref{h et h}) implique :
\begin{eqnarray*}
\mbox{MS}(\nabla (s+t)(x)_{\mbox{Levi}}) & \geq &  \mbox{MS}(\nabla (s+t)(x)_{T_x Y}) - 2 \| \overline{\partial} (s+t)(x) \|
\\  &\geq &  \mbox{MS}(\nabla (s+t)(x)_{T_x Y})  - 2(K+ \varepsilon).
\end{eqnarray*}

La quantit 
$$ \mbox{max} \left\{\left\| (s+t)(x) \right\| ,\; k^{-\frac{1}{2}} \times \mbox{MS}(\nabla (s+t)(x)_{\mbox{Levi}})\right\}  $$
est donc minore par :
\begin{eqnarray*}
 \mbox{max} \left\{\left\| (s+t)(x) \right\| ,\; k^{-\frac{1}{2}} \times \mbox{MS}(\nabla (s+t)(x)_{T_x Y})\right\} - 2(K+ \varepsilon) k^{-\frac{1}{2}}
 & \geq & 2\eta - 2(K+\varepsilon) k^{-\frac{1}{2}}
\\ & \geq & \eta
\end{eqnarray*}
pour $k$ assez grand.

\subsection{Construction de sous-varits de contact}
Soit une varit compacte $Y$ munie d'une forme de contact $\alpha$.
Notons $X$ la symplectise, i.e. le produit $\R \times Y$ muni de la forme symplectique $\omega=d(e^t \alpha)$.
On identifie $Y$  l'hypersurface $\{ 0 \} \times Y$.
Munissons $X$ d'une structure presque complexe $J$, compatible avec $\omega$ et prservant la direction de contact 
en tout point de $Y$, autrement dit la direction de contact sera gale  la direction de Levi.
Comme pr-quantification sur $X$, choisissons le fibr trivial en droites hermitiennes, muni de la connexion suivante :
$$
\nabla =d -i2\pi e^t \alpha.
$$
Comme prcedemment, on suppose donns deux rels $\varepsilon >0$ et $K>0$ et deux entiers $m_{max}  \geq 0$ et $r>0$. 
\\
\\
\indent
Thorme d'Ibort, Mart\'{i}nez-Torres et Presas.

Pour tout entier $k$ assez grand et toute section $(K,2)-$contrle et $(K,1)-$approximativement holomorphe 
$s$ de $rL^k$, il existe une section $(\varepsilon,m_{max} )-$contrle et $(\varepsilon,m_{max} )-$approximativement holomorphe $t$ de $rL^k$ telle que le
lieu des zros de la restriction de $s+t$  $Y$ soit une sous-varit de contact de codimension $2r$ dans $Y$.
\\
\\
\indent
Dmonstration. 
Le thorme (\ref{Levi}) fournit une section approximativement holomorphe $\sigma = s+t$ qui s'annule 
transversalement le long de la direction de Levi (avec une estime). 
Posons $Y_1 = Y\cap \sigma^{-1}(0)$.
Pour $k$ assez grand, on vrifie immdiatement les trois points suivants :

1) L'annulation transversale le long de $Y$ implique que $Y_1$ est une sous-varit lisse 
d'espace tangent $TY \cap \mbox{ker} \nabla \sigma$.

2) La forme $\alpha$ n'est pas nulle sur $TY$ et la restriction de $\nabla \sigma$
 la direction de Levi $TY \cap \mbox{ker }\alpha$ est surjective. Ces deux faits quivalent  la surjectivit, sur $TY$,
de l'application linaire $(\alpha, \nabla \sigma)$, elle-mme quivalente aux deux faits suivants :
$\nabla \sigma$ est surjective sur $TY$ et $\alpha$ n'est pas nulle sur $TY \cap \mbox{ker} \nabla \sigma$
($=TY_1$).

3) L'application $\R-$linaire $\nabla \sigma_{\mbox{Levi}}$ est approximativement $\C-$linaire en ce sens que la
norme de sa partie antilinaire est trs petite compare au module de surjectivit.
On en dduit que le noyau de $\nabla \sigma_{\mbox{Levi}}$ est un sous-espace approximativement $J-$invariant
et, en particulier, symplectique.

Bilan.
Nous avons prouv que $Y_1$ est une sous-varit lisse
sur laquelle la forme $\alpha$ ne s'annule pas et que $(TY_1 \cap \mbox{ker }\alpha, \, d\alpha)$
est un espace vectoriel symplectique donc $Y_1$ est une sous-varit de contact.
\\
\\
\indent
Remarque.
Dans \cite{IbMaPr00}, cette construction que nous avons faite pour $rL^k$ tait faite plus gnralement pour la version twiste $L^k \otimes V$, o $V$ dsigne
un fibr vectoriel hermitien, mais le thorme de Donaldson-Auroux relatif admet bien sr une telle version twiste.

\section{Le thorme d'unicit}\label{chapau}

Ce chapitre a deux objets : premirement dduire le thorme d'unicit d'Auroux (valable dans le cas absolu) du thorme de Donaldson-Auroux relatif
et deuximement donner un contre-exemple qui prouve qu'il n'y a pas unicit dans le cas relatif.

Disons un mot d'une ide sous-jacente dans les preuves des thormes (\ref{Levi}) et (\ref{Auroux}).
En gnral, le fait qu'une famille de sections s'annule transversalement n'implique pas que chaque section de la
famille s'annule transversalement. Pourtant, c'est le cas dans la situation qui nous intresse : si une famille
 $1$ paramtre rel de sections (approximativement) holomorphes s'annule transversalement (avec des estimes),
chaque section s'annulera transversalement.  

\subsection{Sections au-dessus de la sphre de Riemann}\label{sphere}
Notons $L_{\C P^1}$ le fibr hyperplan sur la sphre de Riemann $\C P^1$. C'est une pr-quantification de $\C P^1$ (pour les mtriques
usuelles de la base et de la fibre et la connexion de Chern).

Alors les sections holomorphes de $L_{\C P^1}^k$, pour $k\geq 1$, sont les polynmes 
homognes de degr $k$ en deux variables notes $z_0$ et $z_1$.
Voici quatre exemples de sections de ce fibr.
Les sections polaires :
\begin{eqnarray*}
N_k & = & z_0^k
\\ S_k & = & z_1^k
\end{eqnarray*}
et les sections quatoriales :
\begin{eqnarray*}
E_k &=& 2^{\frac{k}{2}} \times z_0^{\mbox{Ent}\left( \frac{k}{2} \right)} \times z_1^{k-\mbox{Ent}\left( \frac{k}{2} \right)}
\\F_k & = & 2^{\frac{k}{2}} \times z_0^{1+\mbox{Ent}\left( \frac{k}{2} \right)} \times z_1^{k-1-\mbox{Ent}\left( \frac{k}{2} \right)}
\end{eqnarray*}
o Ent dsigne la partie entire.
Alors pour tout entier $m_{max}\geq 0$, il existe un rel $K>0$
indpendant de $k$
tel que ces quatre sections soient $(K,\, m_{max})-$contrles
(cela rsulte d'un calcul que nous omettons).

\subsection{Unicit dans le cas absolu}\label{Auroux}
Thorme d'unicit d'Auroux.

Soient une varit presque khlrienne $X$ compacte et munie d'une pr-quantification $L$, deux rels $K$ et $\eta>0$ et deux entiers $r>0$ et $m_{max}\geq 2$.
Alors il existe un entier $k_0 >0$ et des rel $K' \geq K$ et $\eta' \in ]0,\, \eta[$ tels que, pour tout entier $k\geq k_0$ et toutes sections $s_0$ et $s_1$ de $rL^k$,
si $s_0$ et $s_1$ sont $(K,\, m_{max})-$contrles, $(K,\, m_{max})-$approximativement holomorphes et s'annulent transversalement sur $X$ avec un module 
de transversalit de poids $(1,\, k^{-\frac{1}{2}})$ suprieur  $\eta$,
elles seront isotopes parmi les sections  $(K',\, m_{max})-$contrles, $(K',\, m_{max})-$approximativement holomorphes 
qui s'annulent transversalement avec un module de transversalit de poids $(1,\, k^{-\frac{1}{2}})$ suprieur  $\eta'$.
\\
\\
\indent
Dmonstration.
On munit le produit $\C P^1 \times X$ de la structure presque khlrienne produit 
et de la pr-quantification produit, i.e. le fibr en droites $\overline{L}$ dfini par $\overline{L} = L_{\C P^1} \otimes L$.
On dfinit une section $s$ de $\overline{L}^k$ par :
\begin{eqnarray*}
s & = & N_k \otimes s_0 + S_k \otimes s_1
\end{eqnarray*}  
o $N_k$ et $S_k$ dsignent les sections polaires de $L_{\C P^1}^k$ dfinies au (\ref{sphere}).

Notons $\R P^1$ le cercle mridien form des lments de $\C P^1$ de la forme $[r_0 : r_1]$ avec $r_0$ et $r_1$ rels
et $\R P_+^1$ le demi-cercle mridien form des lments de $\R P^1$ de la mme forme avec $r_0\geq 0$ et $r_1\geq 0$. 
On trivialise $L_{\C P^1}$ au-dessus, de $\R P_+^1$ en identifiant  $1$ la section unitaire $\sqrt{1+ \left( \frac{r_1}{r_0}\right)^2} \times z_0$
(prolonge par la valeur $z_1$ au-dessus du point $[0 : 1]$). Cette trivialisation permet, en un point $(p,\, x) \in \R P_+^1 \times X$,
d'identifier les fibres $rL_x^k$ et $r\overline{L}_{(p,\, x)}^k$. Ainsi, une section de $r\overline{L}^k$ au-dessus de $\R P_+^1 \times X$
pourra tre vue comme une famille  $1$ paramtre de sections de $rL^k$ (au-dessus de $X$).

En appliquant le thorme (\ref{Levi})  l'hypersurface $\R P^1 \times X$ de $\C P^2 \times X$,  la section $s$
et  la prcision $\frac{\eta}{2}$, 
on obtient une section $t$ de $r\overline{L}^k$ qui est 
$(\frac{\eta}{2},\, m_{max})-$contrle, $(\frac{\eta}{2},\, m_{max})-$approximativement holomorphe
et telle que $s+t$ vrifie, pour un certain
rel $\eta''>0$ indpendant de $k$, la minoration suivante :
$$
\mbox{max} \left\{\left\| (s+t)(p,\, x) \right\| ,\; k^{-\frac{1}{2}} \times \mbox{MS}(\nabla (s+t)(p,\, x)_X)\right\}  \geq \eta''
$$
o $\nabla (s+t)(p,\, x)_X$ dsigne la restriction de $\nabla (s+t)(p,\, x)$ au second facteur $T_x X$ de $T_{(p,\, x)} (\C P^1 \times X)$
(en effet, ce facteur $T_x X$ est gal a l'hyperplan complexe de Levi).
On obtient donc une famille  $1$ paramtre de sections de $rL^k$, qui s'annulent transversalement sur $X$
avec un module de poids $(1,\, k^{-\frac{1}{2}})$ suprieur  $\eta''$.
Nous noterons bien sr $s_0+t_0$ et $s_1+t_1$ les deux sections  extrmes de cette famille. Elles sont isotopes parmi les sections de module de 
transversalit pondr suprieur  $\eta''$.

Par ailleurs, isotopons les sections $s_0$ et $s_0 +t_0$ par les sections barycentres $\sigma_{\theta}=s_0 + \theta \times t_0$, avec $\theta \in [0,\, 1]$.
Alors le module de transversalit pondr de $\sigma_{\theta}$ vrifie, en tout point $x\in X$, la minoration suivante :
\begin{eqnarray*}
\\ && \mbox{max} \left\{\left\| \sigma_{\theta}(x) \right\| ,\; k^{-\frac{1}{2}} \times \mbox{MS}(\nabla \sigma_{\theta}(x))\right\} 
\\ & \geq &
\mbox{max} \left\{\left\| s_0(x) \right\|-\theta \times \| t_0 (x) \| ,\; k^{-\frac{1}{2}} \times \left( \mbox{MS}(\nabla s_0(x)) - \theta \times \| \nabla t_0(x) \| \right) \right\} 
\\ && \mbox{(car le module de surjectivit est $1-$lipschitzien)}
\\ & \geq & 
\mbox{max} \left\{\left\| s_0(x) \right\| ,\; k^{-\frac{1}{2}} \times \mbox{MS}(\nabla s_0(x)) \right\} 
- \theta \times \mbox{max} \left\{\| t_0 (x) \| ,\; k^{-\frac{1}{2}} \times \| \nabla t_0(x) \| \right\} 
\\ & \geq & \eta - \theta \times \frac{\eta}{2}
\\ & \geq & \frac{\eta}{2}. 
\end{eqnarray*}
Les sections $s_0$ et $s_0+t_0$ sont donc isotopes parmi les sections de module de 
transversalit pondr suprieur  $\frac{\eta}{2}$. Bien sr, il en va de mme pour $s_1$ et $s_1 +t_1$.

Bilan.
Posons $\eta' = \min \left\{ \frac{\eta}{2}, \eta '' \right\}$.
Les quatre sections $s_0$, $s_0+t_0$, $s_1+t_1$ et $s_1$ 
sont isotopes parmi les sections de $rL^k$ de module de transversalit pondr $\geq \eta'$, ce qui dmontre le thorme.
\\
\\
\indent
Remarque. Dans \cite{Au97}, les hypothses d'holomorphie approximative faites respectivement sur $s_0$ et $s_1$ 
font intervenir deux structures presque complexes $J_0$ et $J_1$ non ncessairement gales (mais compatibles
avec une mme forme symplectique). Cette gnralisation peut tre obtenue en modifiant la dmonstration prcdente.
En effet la sphre de paramtres qui intervenait dans cette dmonstration peut tre utilise
pour interpoler (en un certain sens que nous ne prciserons pas) entre $J_0$ et $J_1$.

\subsection{Un contre-exemple  l'unicit dans le cas relatif}
La partie ${\cal E}$ de la sphre $\C P^1$ d'quation $| z_0 | = | z_1 |$ est appele le cercle quateur.
Les sections quatoriales $E_k$ et $F_k$ dfinies au (\ref{sphere}) tant unitaires au-dessus de ${\cal E}$, leur module de transversalit 
de poids $(1,\, k{^{-\frac{1}{2}}})$ le long de ${\cal E}$ est $\geq 1$ (en fait, il est gal  $1$).

Notons $L_{\cal E}^k$ le fibr $L^k$ restreint au cercle quateur ${\cal E}$. Fixons une trivialisation unitaire de $L_{\cal E}^k$. Dans
cette trivialisation, les restrictions  ${\cal E}$ des deux sections $E_k$ et $F_k$ sont des applications continues de ${\cal E}= S^1$ vers $S^1$.
Comme leurs degrs topologiques diffrent d'une unit, ces deux applications ne sont  pas isotopes parmi les 
applications de $S^1$ vers $\C$ qui ne s'annulent pas.

Remarquons que la dimension relle de ${\cal E}$ (gale  $1$) est strictement infrieure au rang rel du fibr $L_{\cal E}^k$
(gal  $2$) et que par consquent, pour une section de ce fibr, {\it s'annuler transversalement} signifie {\it ne pas s'annuler}.
Les sections $E_k$ et $F_k$ ne sont donc pas isotopes parmi les sections de $L^k$ qui s'annulent transversalement le long de ${\cal E}$. 

Ce contre-exemple prouve que le thorme d'unicit d'Auroux, valable dans le cas absolu, n'admet pas d'analogue dans le cas relatif.

\section{Variations de Vitushkin}\label{chapvv}

Ce chapitre dcrit la thorie des variations de Vitushkin et l'applique  l'tude de la complexit des hypersurfaces algbriques
(il est inspir de \cite{CoYo04}). 

\subsection{Partie espace}
Dans un espace euclidien, une partie finie sera dite $\varepsilon -${\it espace} si deux quelconques de ses points vrifient : $d(x,\, y) \geq \varepsilon $.

\subsection{Partie espace dans une hypersurface}
Rappelons sans dmonstration un rsultat lmentaire de gomtrie des hypersurfaces :
\\
\\
\indent
Proposition. Soit une hypersurface $A$  compacte  bord dans un espace euclidien de dimension $n$. 
Alors il existe une constante $C$ telle que pour tout $\varepsilon \in ]0, 1[$
et toute partie finie  $\varepsilon -$espace $F$ de $A$ 
le cardinal de $F$ soit major par $C \varepsilon^{1-n}$.

\subsection{Partie espace dans une hypersurface algbrique}\label{comp hyp}
On veut comprendre comment $C$ dpend de $A$, dans le cas d'une hypersurface algbrique.
La rponse est donne par le thorme suivant :
\\
\\
\indent
Thorme. 
Soit un espace euclidien $E$ de dimension $n$.
Notons $B_E$ la boule unit ferme de $E$.
Alors il existe un polynme $P : \R \mapsto \R$
tel que pour toute hypersurface algbrique $A \subset E$, tout rel $\varepsilon  \in ]0, 1[$,
et toute partie finie $\varepsilon -$espace $F$ de $A \cap B_E$, le cardinal
de $F$ soit major par $P( \mbox{deg }A) \varepsilon^{1-n}$.

Cette majoration est notamment vraie pour une $\varepsilon -$discrtisation, c'est  dire une 
partie $\varepsilon -$espace maximale.
\\
\\
\indent
L'outil qui nous permettra de dmontrer ce thorme s'appelle les variations de Vitushkin.
Nous allons les dfinir au (\ref{vdv}) puis en donner quelques proprits et en dduire le thorme. 

\subsection{Mesure sur une grassmannienne linaire}
Soit $E$ un espace euclidien.
Notons $\mbox{Gr}(k,\, E)$ la grassmanienne des sous-espaces vectoriels de dimension $k$.
Il existe sur $\mbox{Gr}(k,\, E)$ une mesure invariante par isomtries linaires, unique  multiplication par un rel prs.
Notons $\mbox{Mes}(k,\, E)$ cette mesure, normalise par la condition :
$$
\int_{\mbox{\small Gr}(k,\, E)} 1 = 1.
$$
Par convention, on supposera toujours $\mbox{Gr}(k,\, E)$ munie de cette mesure.

\subsection{Mesure sur une grassmannienne affine}
Notons $\mbox{AffGr}(k,\, E)$ la grassmanienne des sous-espaces affines de $E$ de dimension $k$.
Si $F$ dsigne un sous-espace vectoriel de $E$ et $p$ un point de l'orthogonal $F^{\perp}$, on notera $p+F$
le sous-espace affine de $E$ parallle  $F$ passant par $p$.
Par  convention, on supposera toujours $\mbox{AffGr}(k,\, E)$ munie de la mesure $\mbox{AffMes}(k,\, E)$ 
dfinie par la relation :
$$
\int_{F \in \mbox{\small AffGr}(k,\, E)} \phi (F)\, dF = \int_{F \in \mbox{\small Gr}(k,\, E)} \left( \int_{p\in F^{\perp}} \phi(p+F) \, dp \right) \, dF.
$$
Cette relation est valable pour toute fonction mesurable $\phi :  \mbox{AffGr}(k,\, E) \rightarrow [0,\, + \infty]$.
\\
\\
\indent
Remarque.
Cette mesure $\mbox{AffMes}(k,\, E)$ est l'unique mesure sur $\mbox{AffGr}(k,\, E)$ invariante par isomtries affines,
 multiplication par un rel prs.

\subsection{La formule de Cauchy-Crofton}

Thorme.
Soit un espace euclidien $E$ de dimension $n$.
Alors toute sous-varit $X$ de dimension $d$ de $E$ vrifie :
$$
\int_{F \in \mbox{\small AffGr}(n-d,\, E)} \mbox{Card} (X \cap F) \, dF = c_{d,\, n} \mbox{Vol}_d X
$$
o $c_{d,\, n}$ dsigne un rel strictement positif qui ne dpend que des dimensions $d$ et $n$,
o Card$(X \cap F)$ dsigne le cardinal de l'intersection $X \cap F$ et o $\mbox{Vol}_d X$ dsigne le 
volume ($d-$dimensionnel) de la sous-varit $X$.
\\
\\
\indent
Dmonstration.

Pour tout sous-espace vectoriel $F$ de $E$ de dimension $d$, notons $p_X^{F}$ la projection orthogonale de $X$ vers $F$.
En un point $x$ de $X$, sa diffrentielle, note $p_{\, T_x X}^{F}$, est une application linaire entre l'espace tangent $T_x X$ et le sous-espace $F$
qui sont deux espaces vectoriels de mme dimension $d$. 
Notons $\Lambda ^d \,  p_{\, T_x X}^{F}$ la dernire puissance extrieure de l'application linaire $p_{\, T_x X}^{F}$,
c'est--dire l'application linaire induite entre les deux droites vectorielles euclidiennes
$\Lambda ^d\, T_x X$ et $\Lambda ^d F$.
Nous allons calculer de deux faons diffrentes l'intgrale double suivante :
$$
I = \int_{(x,\, F) \in X\times \mbox{\small Gr}(d,\, E)} \left\| \Lambda ^d \,  p_{\, T_x X}^{F} \right\| \, dx\, dF.
$$

Pour tout $F \in \mbox{Gr} (d,\, E)$, appliquons la formule de changement de variables :
$$
\int_{x \in X}  \left\| \Lambda ^d \,  p_{\, T_x X}^F \right\| \, dx 
= \int_{y \in F} \mbox{Card } \left(p_X^{F}\right)^{-1} \{ y \}\, dy
= \int_{y \in F} \mbox{Card } (X \cap (y+ F^{\perp})) \, dy.
$$

Intgrons sur $ \mbox{Gr} (d,\, E)$ :
\begin{eqnarray*}
I & = & \int_{F \in \mbox{\small Gr}(d,\, E)}\left(  \int_{y \in F} \mbox{Card } (X \cap (y+ F^{\perp})) \, dy \right) \, dF
\\ & = & \int_{F \in \mbox{\small Gr}(n-d,\, E)}\left(  \int_{y \in F^{\perp}} \mbox{Card } (X \cap (y+ F)) \, dy \right) \, dF
\\ & = & \int_{F \in \mbox{\small AffGr}(n-d,\, E)} \mbox{Card} (X \cap F) \, dF .
\end{eqnarray*}

Pour $F' \in \mbox{Gr} (d,\, E)$, posons :
$$
c_{d,\, n} = \int_{F \in \mbox{\small Gr}(d,\, E)} \left\| \Lambda ^d \,  p_{\,F'}^F \right\| \, dF.
$$
o $ p_{\,F'}^F $ dsigne la projection orthogonale de $F'$ vers $F$.
Ce rel $c_{d,\, n}$ ne dpend pas de $F'$ car le groupe des isomtries linaires opre
transitivement sur la grassmannienne $\mbox{Gr} (d,\, E)$.
En prenant $F' = T_x X$ et en intgrant sur $X$, on obtient :
$$
I = \int_{x\in X} c_{d,\, n} \, dx = c_{d,\, n} \mbox{Vol}_d X.
$$

Il suffit alors de comparer les deux expressions de $I$ obtenues pour achever la dmonstration du thorme.

\subsection{Les variations de Vitushkin dans le cas absolu}\label{vdv}
Soit $E$ un espace euclidien de dimension $n$.
Les variations de Vitushkin 
d'une partie compacte $A$ de $E$ sont des sortes de mesures $d$ dimensionnelles de cet ensemble
pour $d$ compris entre $0$ et $n$.
Plus prcisment :
\\
\\
\indent
Dfinition.
La $d-${\it me variation de Vitushkin} de $A$ est l'intgrale sur la grassmanienne 
des sous-espaces affines de dimension $n-d$ de $E$ du nombre de composantes connexes
de $A \cap F$. On la note $V^{E}_d(A)$ ou plus simplement $V_d(A)$.
\\
\\
\indent
Exemple 1 :
Le nombre $V_0(A)$ est le nombre de composantes connexes de $A$.

Exemple 2 :
Le nombre $V_n(A)$ est le volume de $A$.

Exemple 3 :
D'aprs la formule de Cauchy-Crofton, toute sous-varit $A$ de dimension $d$
vrifie :
$$
V_d(A) = c_{d,\, n}\mbox{Vol}_d(A)
$$
o $\mbox{Vol}_d(A)$ dsigne le volume $d-$dimensionnel de $A$
et o $c_{d,\, n}>0$ dsigne un rel ne dpendant que des dimensions $d$ et $n$.

\subsection{Les variations de Vitushkin dans le cas relatif}

Nous utiliserons plutt la version relative :
\\
\\
\indent
Dfinition.
Soient, dans $E$, une partie compacte $A$ et une partie ferme $B$. Notons $d$ un entier compris entre $0$ et $n$.
La $d-${\it me variation de Vitushkin} du couple $(A,\, B)$ est l'intgrale sur la grassmanienne 
des sous-espaces affines de dimension $n-d$ de $\R^n$ du nombre de composantes connexes
de $A \cap F$ disjointes de $B$. On la note $V_d^{E}(A,\, B)$ ou plus simplement $V_d(A,\, B)$.
\\
\\
\indent
Exemple 1 :
Le nombre $V_0(A,\, B)$ est le nombre de composantes connexes de $A$ disjointes de $B$.

Exemple 2 :
Le nombre $V_n(A,\, B)$ est le volume de $A \backslash B = \{ a\in A | a\notin B\}$.

\subsection{Un lemme technique}\label{lemme tech}
Soient un espace topologique $X$, une partie $Y\subset X$ et un point $y\in Y$.
Notons $C_X(y)$ et $C_Y(y)$ les composantes connexes respectives de $y$ dans $X$ et dans $Y$.
En gnral, elles vrifient l'inclusion :
$$
C_Y(y) \subset C_X(y).
$$
L'autre inclusion est vrifie sous certaines hypothses.
\\
\\
\indent
Lemme.
On suppose que la partie $Y$ est compacte et que la composante connexe $C_Y(y)$
est disjointe du bord de $Y$.
Alors les composantes $C_X(y)$ et $C_Y(y)$ concident.
\\
\\
\indent
Dmonstration.
Comme $Y$ est compacte, la composante connexe $C_Y(y)$ est l'intersection des parties  la fois fermes et ouvertes dans $Y$
qui contiennent le point $y$. Ecrivons :
$$
C_Y(y) = \bigcap_i Z_i.
$$
Le bord de $Y$ tant compact et disjoint de $C_Y(y)$, il est disjoint d'une intersection finie $Z = Z_{i_1} \cap \dots \cap Z_{i_k}$.
Cette intersection $Z$ est ferme dans le ferm $Y$ donc dans $X$.
De mme $Z$ est ouverte dans l'intrieur de $Y$ donc dans $X$.
Enfin $Z$ contient le point $y$. De ces trois proprits, on dduit que la composante $C_X(y)$
est incluse dans $Z$ et, par transitivit de inclusion, dans $Y$. Comme $C_X(y)$ est connexe, on en dduit l'inclusion :
$$
C_X(y) \subset C_Y(y)
$$
et la proposition en dcoule car, comme nous l'avons dit, l'autre inclusion est vidente.

\subsection{Additivit des variations de Vitushkin}\label{additif}

Proposition.
Soient des boules fermes $B_i \subset E$ en nombre fini, d'intrieurs disjoints et une partie compacte $A$ de $E$.
Notons $S_i$ les sphres bordant les $B_i$, notons : $\mbox{int}B_i = B_i \backslash S_i$ l'intrieur de $B_i$ 
et posons : $$B = E \backslash \bigcup_i \mbox{int}B_i.$$

Alors on a:
$$
V_d(A,\, B) = \sum_i V_d(A\cap B_i,\, S_i) 
$$
pour tout entier $d$ compris entre $0$ et $n$.
\\
\\
\indent
Dmonstration.
Il suffit de traiter le cas $d=0$. En effet, le cas gnral s'y ramne grce  la formule suivante :
$$
V_d^{E}(X,\, Y) = \int_{F \in \mbox{\small AffGr}(n-d,\, E)} V_0^{F}(X\cap F,\, Y \cap F)\, dF 
$$
qui est une simple reformulation de la dfinition de $V_d^{E}(X,\, Y)$.

Toute composante connexe de $A$ disjointe de $B$ est incluse dans la runion des boules ouvertes 
$\mbox{int}B_i$ et donc, par connexit, dans l'une de ces boules. D'aprs le lemme (\ref{lemme tech}),
les composantes connexes de $A$ incluses dans la boule ouverte $\mbox{int}B_i$ sont exactement les composantes 
connexes de $A \cap B_i$ incluses dans $\mbox{int}B_i$.
En les comptant, on obtient l'identit recherche :
$$
V_0(A,\, B) = \sum_i V_0(A\cap B_i,\, S_i). 
$$

\subsection{Majoration des variations de Vitushkin dans le cas d'une hypersurface algbrique}\label{majore}

Proposition.
Notons $B_E$ la boule unit ferme de l'espace euclidien $E$.
Soit un entier $d$ compris entre $0$ et $n$. Il existe un polynme $P_d : \R \rightarrow \R$ vrifiant
pour toute hypersurface algbrique $A$
et toute partie ferme $B$ de $E$, l'ingalit :
$$
V_d(A\cap B_E,\, B) \leq P_d(\mbox{deg }A).
$$
\\
\indent
Dmonstration. Soit $F$ un sous-espace affine. Notons $f(X)$ le nombre de composantes connexes 
de $A \cap B_E \cap F$. Il est major par une fonction polynomiale de deg $A$.
Le support de $f$ est inclus dans l'ensemble des sous-espaces qui intersectent la
boule $B_E$ et cet ensemble est un domaine de mesure fini. Donc l'intgrale de $f$ est 
majore par une fonction polynomiale de deg $A$. 

\subsection{Hyperplans et intersections}\label{hyp int}
Proposition.
Soient $A$ une partie compacte et $B$ une partie ferme de l'espace euclidien $E$.
On suppose $A$ non incluse dans $B$ et on  suppose que le nombre
$V_0(A,\, B)$ est nul.
Notons $\varepsilon$ et $a$ un rel strictement positif et un point de $A$ vrifiant :
$$
\forall b\in B \;\;\;\;\;\;\;\; d(a,\, b)\geq \varepsilon.
$$ 

Notons ${\cal H}$ la partie de $\mbox{Gr}(n-1,\, E)$ contenant les hyperplans affines $H$ pour lesquels l'intersection  $A \cap H$ contient un point $a_H$
vrifiant :
$$
d(a_H,\, a) \leq \varepsilon.
$$

Alors la mesure de cet ensemble ${\cal H}$ vrifie la minoration suivante :
$$\mbox{mes}({\cal H}) \geq c_{1,n}\varepsilon,$$ 
la constante $c_{1,\, n}$ tant celle de la formule de Cauchy-Crofton.
\\
\\
\indent
Dmonstration.
Notons $A'$ l'ensemble des points $x$ de $A$ qui vrifient $d(a,\, x) \leq \varepsilon$
et notons $S$ la sphre des points $x$ de $E$ qui vrifient $d(a,\, x) = \varepsilon$ .
L'hypothse : $$V_0(A,\, B)=0$$ implique :
$$V_0(A',\, S)=0$$
(cette implication est une consquence du lemme (\ref{lemme tech})).
Il existe donc  un point $b\in S$ appartenant  la mme composante connexe
de $A'$ que le point $a$. 
Tout hyperplan sparant $a$ et $b$ rencontre $A'$ et appartient donc  ${\cal H}$. 
Notons ${\cal K}$ l'ensemble des hyperplans sparant $a$ et $b$.
Appliquons la formule de Cauchy-Crofton au segment $[a\; b]$  :
$$\mbox{mes}({\cal K}) = c_{1,\, n}d(a,\, b) = c_{1,\, n}\varepsilon.$$
Alors, l'inclusion de ${\cal K}$ dans ${\cal H}$ permet de conclure.

\subsection{Changement d'espace euclidien}\label{baisse}
Proposition.
Soient un espace euclidien $E$ de dimension $n$, un sous-espace affine $F$ de dimension $k$ et un entier naturel $d\leq k$.
Notons $U$ l'ouvert de $\mbox{AffGr}(n-d,\, E)$ form des sous-espaces affines de direction transversale  celle de $F$.
Alors toute fonction mesurable $\phi :  \mbox{AffGr}(k-d,\, F) \rightarrow [0,\, + \infty]$ vrifie la relation:
$$
\int_{X\in  \mbox{AffGr}(k-d,\, F)} \phi(X)\, dX =c_{d,\, k,\, n} \int_U \phi (X\cap F)\, dX
$$
o $c_{d,\, k,\, n}$ dsigne une constante strictement positive ne dpendant que des dimensions.
\\
\\
\indent
Dmonstration.
Il suffit de remarquer que les deux membres de cette identit s'interprtent comme les
intgrales de la fonction $\phi$ sur $\mbox{AffGr}(k-d,\, F)$ par rapport  des mesures 
invariantes par les isomtries affines de $F$ et que deux telles mesures sont gales, 
 multiplication prs par un rel.
\\
\\
\indent
Corollaire. Soient $A$ une partie compacte et $B$ une partie ferme du sous-espace  affine $F$.
Alors les $d-$mes variations de Vitushkin de la paire $(A,\, B)$ dans $E$ et dans $F$
vrifient la relation suivante :
$$
V_d^{F}(A,\, B) = c_{d,\, k,\, n} \, V_d^{E}(A,\, B).
$$
\\
\\
\indent
Dmonstration.
Pour tout sous-espace affine $X$ de $E$, notons $\phi(X)$
le nombre de composantes connexes de $A\cap X$ disjointes de $B$. 
\begin{eqnarray*}
V_d^{F}(A,\, B) &=&
\int_{X\in  \mbox{\small AffGr}(k-d,\, F)} \phi(X)\, dX 
\\ 
& = & c_{d,\, k,\, n} \int_{X\in U} \phi (X\cap F)\, dX
\\
& = &
c_{d,\, k,\, n} \int_{ X\in \mbox{\small AffGr}(n-d,\, E)} \phi (X\cap F)\, dX
\\ & = &
c_{d,\, k,\, n} \, V_d^{E}(A,\, B)
\end{eqnarray*}

La deuxime identit est une application de la proposition et 
la troisime identit dcoule du fait que la mesure du 
 complmentaire de $U$ dans $ \mbox{AffGr}(n-d,\, E)$ est nulle.
 
 \subsection{Une application du thorme de Fubini}\label{fubi}

Proposition.
Soient deux entiers naturels $k_1$ et $k_2$ avec $k_1 \leq k_2 \leq \mbox{dim }E$.
Alors toute fonction mesurable $\phi :  \mbox{AffGr}(k_1,\, E) \rightarrow [0,\, + \infty]$ vrifie 
la relation suivante :
$$
\int_{F_1 \in \mbox{AffGr}(k_1,\, E)} \phi (F_1)\, dF_1 = \int_{F_2 \in \mbox{AffGr}(k_2,\, E)} \left( \int_{F_1 \in \mbox{AffGr}(k_1,\, F_2)} \phi (F_1)\, dF_1\right)\, dF_2.
$$

Dmonstration.
Commenons par indiquer une formule analogue pour les grassmanniennes linaires.
Toute fonction mesurable $\phi :  \mbox{Gr}(k_1,\, E) \rightarrow [0,\, + \infty]$ vrifie :
$$
\int_{F_1 \in \mbox{Gr}(k_1,\, E)} \phi (F_1)\, dF_1 = \int_{F_2 \in \mbox{Gr}(k_2,\, E)} \left( \int_{F_1 \in \mbox{Gr}(k_1,\, F_2)} \phi (F_1)\, dF_1\right)\, dF_2.
$$
Pour le voir, il suffit de remarquer que les deux membres de cette identit s'interprtent comme les
intgrales de la fonction $\phi$ sur $\mbox{Gr}(k_1,\, E)$ par rapport  des mesures et que ces mesures sont gales par unicit
sur $\mbox{Gr}(k_1,\, E)$ de la mesure normalise invariante par isomtries linaires.

Passons au cas affine. 
Une fonction mesurable $\phi :  \mbox{AffGr}(k_1,\, E) \rightarrow [0,\, + \infty]$ vrifie :
\begin{eqnarray*}
&&\int_{F_1 \in \mbox{AffGr}(k_1,\, E)} \phi (F_1)\, dF_1 
\\
& = & 
\int_{F_1 \in \mbox{Gr}(k_1,\, E)} \left( \int_{p_1 \in F_1^{\perp}} \phi (p_1 + F_1)\, dp_1 \right) \, dF_1 
\\
& = & 
\int_{F_2 \in \mbox{Gr}(k_2,\, E)} \left( \int_{F_1 \in \mbox{Gr}(k_1,\, F_2)} \left( \int_{p_1 \in F_1^{\perp}} \phi (p_1 + F_1)\, dp_1 \right) \, dF_1 \right) \, dF_2 
\\
& = & 
\int_{F_2 \in \mbox{Gr}(k_2,\, E)} \left( \int_{F_1 \in \mbox{Gr}(k_1,\, F_2)} \left( \int_{(p_1,\, p_2) \in (F_1^{\perp} \cap F_2) \oplus F_2^{\perp}} 
\phi (p_1+p_2 + F_1)\, dp_1 \, dp_2\right) \, dF_1 \right) \, dF_2 
\\
& = & 
\int_{F_2 \in \mbox{Gr}(k_2,\, E)} \left( \int_{p_2 \in F_2^{\perp}} \left( \int_{F_1 \in \mbox{AffGr}(k_1,\, F_2)} 
\phi (p_2 + F_1)\, dF_1\right) \, dp_2 \right) \, dF_2 
\\
& = & 
\int_{F_2 \in \mbox{AffGr}(k_2,\, E)} \left( \int_{F_1 \in \mbox{AffGr}(k_1,\, F_2)} 
\phi (F_1)\, dF_1\right) \, dF_2 .
\end{eqnarray*}
La dernire identit dcoule du fait que la translation de vecteur $p_2$ induit un isomorphisme d'espaces mesurs
entre $\mbox{AffGr}(k_1,\, F_2)$ et $\mbox{AffGr}(k_1,\, p_2 + F_2)$.
La proposition est dmontre.
 \\
 \\
 \indent
 Corollaire.
 Soient $A$ une partie compacte et $B$ une partie ferme de $E$ 
 et soient $k$ et $d$ deux entiers naturels vrifiant $0\leq n-k \leq d \leq n$.
 Alors :
 $$
 V_{d}^{E}(A,\, B) = \int_{F \in \mbox{AffGr}(k,\, E)} V_{d+k-n}^{F}(A \cap F,\, B \cap F)\, dF.
 $$
 
 Dmonstration.
 Il suffit d'appliquer le thorme pour $k_1 = n - d$, pour $k_2 = k$ et pour la fonction $\phi$ qui associe
  un sous-espace affine $X$ de $E$ de dimension $n-d$ le nombre de composantes connexes de $A\cap X$ disjointes de $B$.

\subsection{Minoration de la somme variations de Vitushkin}\label{la somme bis}
Enonons le thorme principal de cette thorie.
\\
\\
\indent
Thorme.
Soient $A$ une partie compacte et $B$ une partie ferme de l'espace euclidien $E$.
On suppose $A$ non incluse dans $B$. Notons $\varepsilon$ un rel strictement positif vrifiant :
$$
\exists a\in A \;\;  \forall b\in B \;\;\;\;\;\;\;\; d(a,\, b)\geq \varepsilon.
$$ 

Alors la minoration suivante est vrifie :
$$
\alpha_n \leq \sum_{d=0}^{n} \, \varepsilon^{-d} \; V_d^{E}(A,\, B)
$$
o $\alpha_n$ dsigne un rel strictement positif ne dpendant que de la dimension $n$ de l'espace $E$.
\\
\\
\indent
Dmonstration par rcurrence.

En dimension $0$, c'est clair. Supposons le thorme dmontr en dimension $n - 1$.

Par hypothse, nous pouvons choisir un point $a \in A$ vrifiant :
$$
\forall b\in B \;\;\;\;\;\;\;\; d(a,\, b)\geq \varepsilon.
$$ 

Nous allons distinguer deux cas selon que $V_0^{E}(A,\, B)$ est nul ou non.
Il est clair que, si nous obtenons dans chacun des deux cas une ingalit du type recherch
(avec des rels strictement positifs respectifs $\alpha_n'$ et $\alpha_n''$), le thorme sera dmontr
(il suffira de poser $\alpha_n = \mbox{min} \{ \alpha_n',\, \alpha_n''\}$).

Le cas $V_0^{E}(A,\, B)\neq 0$ est trs facile. Comme $V_0^{E}(A,\, B)$ est un entier strictement positif, on peut crire :
$$
1 \leq V_0^{E}(A,\, B) \leq \sum_{d=0}^{n} \, \varepsilon^{-d} \; V_d^{E}(A,\, B)
$$
et le rel $\alpha_n' = 1$ convient.

Dsormais, on suppose $V_0^{E}(A,\, B)$ nul.
Comme prcdemment, nous noterons
${\cal H}$ l'ensemble des hyperplans affines $H$ pour lesquels l'intersection  $A \cap H$ contient un point $a_H$
vrifiant :
$$
d(a_H,\, a) \leq \frac{\varepsilon}{2}.
$$

Par ingalit triangulaire  :
$$
\forall b\in B \;\;\;\;\;\;\;\; d(a_H,\, b)\geq d(a,\, b) - d(a_H,\, a) \geq \varepsilon - \frac{\varepsilon}{2} \geq \frac{\varepsilon}{2} .
$$

Appliquons l'hypothse de rcurrence (dans l'hyperplan $H$) et le corollaire de la proposition (\ref{baisse}) :
\begin{eqnarray*}
\alpha_{n-1} &\leq & \sum_{d=0}^{n-1} \, \left(\frac{\varepsilon}{2} \right)^{-d} \; V_d^{H}(A \cap H,\, B\cap H)
\\ & = & \sum_{d=0}^{n-1} \,c_{d,\, n-1,\, n} \left(\frac{\varepsilon}{2} \right)^{-d} \; V_d^{E}(A \cap H,\, B\cap H).
\end{eqnarray*}

En intgrant sur $ {\cal H}$ et en appliquant la proposition (\ref{hyp int}) et le corollaire de la proposition (\ref{fubi}), on obtient :
\begin{eqnarray*}
\frac{\alpha_{n-1} c_{1,\, n} \varepsilon}{2} &\leq & \sum_{d=0}^{n-1} \,c_{d,\, n-1,\, n} \left(\frac{\varepsilon}{2} \right)^{-d} \; \int_{H \in  {\cal H}}V_d^{E}(A \cap H ,\, B\cap H)\, dH
\\ &\leq & \sum_{d=0}^{n-1} \,c_{d,\, n-1,\, n} \left(\frac{\varepsilon}{2} \right)^{-d} \;
V_{d+1}^{E} (A,\, B).
\end{eqnarray*}

On obtient bien une ingalit du type :
$$
\alpha_n'' \leq \sum_{d=1}^{n} \, \varepsilon^{-d} \; V_d^{E}(A ,\, B)
$$
avec $\alpha_n'' >0$. Comme nous l'avons dit, cela suffit  dmontrer le thorme.

\subsection{Le cas d'un ensemble de mesure nulle}\label{cas mes nulle}
Proposition.
Soient $A$ une partie ferme de l'espace euclidient $E$ 
et $B_r$ une boule ferme de $E$, de rayon $r$ compris entre $0$ et $1$ et
centre en un point de $A$. Notons $S_r$ la sphre bordant $B_r$.
On suppose $A$ de volume nul.
Alors la minoration suivante est vrifie :
$$
\alpha_n r^{n-1} \leq \sum_{d=0}^{n-1} V_d(A \cap B_r,\, S_r)
$$
o $\alpha_n$ dsigne un rel strictement positif ne dpendant que de la dimension $n$ de l'espace $E$
\\
\\
\indent
Dmonstration.
Le thorme (\ref{la somme bis}) donne :
$$ \alpha_n  \leq  \sum_{d=0}^{n} r^{-d} V_d(A \cap B_r,\, S_r). $$
Le dernier terme est nul car le volume de $A$ est suppos nul.
Pour les autres termes, l'ingalit $r^{n-1} \leq r^d$ permet de conclure.

\subsection{Dmonstration du thorme (\ref{comp hyp})}
Pour tout point $x$ de $F$, notons $B_x$ le boule ferme de rayon $\frac{\varepsilon}{2}$ centre en $x$
et notons $S_x$ la sphre bordant $B_x$.
Comme la partie $F$ est $\varepsilon-$espace, les intrieurs des boules $B_x$ seront disjoints.
Partitionnons l'ensemble $F$ en deux sous-ensembles $F'$ et $F''$ dfinis respectivemement
par les ingalits $\| x \| \leq 1- \frac{\varepsilon}{2} $ et $ 1- \frac{\varepsilon}{2} < \| x \| $.
 
La proposition (\ref{additif}) implique :
$$
V_d(A\cap B_E,\, B)= \sum_{x\in F'} V_d\left(A \cap B_x,\, S_x\right)
$$
o $B$ dsigne le complmentaire de la runion des intrieurs des $B_x$, $x\in F'$.
La proposition (\ref{majore}) donne :
$$
V_d(A\cap B_E,\, B) \leq P_d(\mbox{deg }A).
$$
Et la proposition (\ref{cas mes nulle}) donne :
$$
\alpha_n \left(\frac{\varepsilon}{2}\right)^{n-1} \leq \sum_{d=0}^{n-1} V_d (A \cap B_x,\, S_x).
$$
En assemblant ces trois rsultats, on obtient :
$$
\alpha_n \left(\frac{\varepsilon}{2}\right)^{n-1} \mbox{Card }F' 
\leq \sum_{d=0}^{n-1} \sum_{x\in F'}V_d (A \cap B_x,\, S_x)
=   \sum_{d=0}^{n-1} V_d(A\cap B_E,\, B) \leq  \sum_{d=0}^{n-1} P_d(\mbox{deg }A).
$$
donc Card $F'$ est major par $\varepsilon^{1-n}P(\mbox{deg }A)$ o $P$ dsigne un polynme. 

Par ailleurs le cardinal de $F''$ est facile  majorer. En effet $F''$ est une partie $\varepsilon-$espace
et les boules $B_x$ centres en les points $x$ de $F''$ sont incluses dans la couronne dfinie par l'encadremement
$1-\varepsilon \leq \|x\| \leq 1+ \frac{\varepsilon}{2}$. Comme le volume de cette couronne est major par une fonction linaire de $\varepsilon$
et que les volumes des boules sont en $\varepsilon^n$, le cardinal de $F''$ est major par une fonction linaire de $\varepsilon^{1-n}$
indpendante de $A$.

Le cardinal de $F=F' \cup F''$ est donc major et
le thorme (\ref{comp hyp}) est dmontr.

\section{Hypersurface de Bertini proche des images des points $\varepsilon-$critiques}\label{chapbert}

\subsection{Deux noncs}

\subsubsection{L'nonc du thorme}\label{lezen}

Soit une application polynomiale $F: \R^m \rightarrow \R^n$ et un $\varepsilon >0$.
Un point de $\R^m$ est dit $\varepsilon-$critque si le module de surjectivit de la diffrentielle de $F$
en ce point est major par $\varepsilon$.

Le but de ce chapitre est de dmontrer le rsultat suivant :
\\
\\
\indent
{\bf Thorme.}
Il existe une hypersurface $H \in \R^n$, dfinie par un polynme non nul $h :  \R^n \rightarrow \R$
de degr $\leq P(D)$ telle que pour tout $X\in \R^m$, $\varepsilon-$critique de norme $\leq 1$,
il existe un point de $H$ qui est $\varepsilon P(D)-$proche de $F(X)$.
\\
 \\
 \indent
 Convention. Dans tout ce chapitre, on note $D$ un majorant du degr des donnes 
 et $P(D)$ un polynme en $D$ (les dimensions tant fixes), qui peut varier d'un nonc  l'autre.
\\
\\
\indent
L'article \cite{Do99} contient un rsultat analogue (pour les applications polynomiales complexes) et en donne une dmonstration
que nous allons suivre.  

\subsubsection{L'nonc d'une gnralisation}\label{lezen 2}
Afin de dmontrer le thorme (\ref{lezen}) par rcurrence sur la dimension $n$ de l'espace d'arrive de $F$, ajoutons-lui des paramtres et des contraintes.

Dsormais, les donnes seront donc deux applications polynomiales $F: \R^m \times \R^p \rightarrow \R^n$
et $G: \R^m \times \R^p \rightarrow \R^q$ et un rel $\varepsilon >0$.
Si $n$ est gal  $1$, la fonction $F$ sera note $f$.
Les indtermines de $F(X,\, T)$ (qui sont aussi celles de $G(X,\, T)$) se sparent en deux familles
que nous appelerons les variables $X$ et les paramtres $T$.
Voici le thorme gnralis dont la dmonstration par rcurrence occupera toute la suite de ce chapitre :
\\
\\
\indent
{\bf Thorme gnralis.}
Il existe une hypersurface $H \in \R^n \times \R^p$, dfinie par un polynme non nul $h :  \R^n \times \R^p \rightarrow \R$
de degr $\leq P(D)$ telle que pour tout paramtre $T\in \R^p$, si $X$ dsigne un point  de $\R^m$ de norme $\leq 1$
en lequel :
\\
\\
\indent
1) $G(.,\, T)$ s'annule transversalement, 

2) le module de surjectivit du gradient de la restriction de $F(.,\, T)$ au lieu d'annulation de $G(.,\, T)$ est major par $\varepsilon$,
\\
\\
\indent
alors il existe un point $Y\in \R^n$, $\varepsilon P(D)-$proche de $F(X,\, T)$ tel que $(Y,\, T)$ appartienne  $H$.

\subsubsection{Gnrique implique gnral}\label{gen gen}
Affirmation. Il suffit de pouver le thorme (\ref{lezen 2}) pour presque tous $F$, $G$, $\varepsilon$ et $T$ (pour un $n$ donn).
\\
\\
\indent
Remarque. {\it Presque tout $F$ (ou $G$)} signifie {\it presque tout $F$ (resp. $G$) parmi les applications polynomiales de degr $\leq D$}.
\\
\\
\indent
Dmonstration de l'affirmation.
Soient des donnes $F$, $G$, $X$, $T$ et $\varepsilon$ satisfaisant les hypothses du thorme.
Donnons-nous une suite $(F_n,\, G_n,\, T_n)$ qui converge vers $(F,\, G,\, T)$
et une suite $(\varepsilon_n)$ qui converge vers $2\varepsilon$
telles le thorme soit vrifi pour les donnes $(F_n,\, G_n,\, \varepsilon_n,\, T_n)$
sur une boule concentrique de rayon $2$, disons.
Alors, le thorme fournit une hypersurface $H_n$ dfinie par un polynme $h_n$ non nul de degr major par $P(D)$. 
Comme le degr est major, on peut (quitte  extraire) supposer que la suite $(h_n)$
converge projectivement vers un polynme $h$ non nul de degr major par $P(D)$. 
Notons $H$ l'hypersurface dfinie par $h$.

Une annulation transversale est ncessairement stable. Il existera donc une suite $(X_n)$
convergeant vers $X$ telle que, pour tout $n$ assez grand, l'application $G_n$ s'annule au point $(X_n,\, T_n)$.
Notons $Y_n$ un point de $\R^n$, $\varepsilon_n P(D)-$proche de $F_n(X_n,\, T_n)$ tel que $(Y_n,\, T_n)$ appartienne  $H_n$.
La suite $(Y_n)$ est borne et donc, quitte  extraire, elle converge. Sa limite $Y$ est 
$2\varepsilon P(D)-$proche de $F(X,\, T)$ et le point $(Y,\, T)$ appartient  $H$.

\subsection{Le cas $n=1$.}
Pour initialiser la rcurrence, on suppose $n$ gal  $1$.

\subsubsection{La sous-varit $A$}
Un paramtre $T$ sera dit $G-$rgulier si l'application $G(.,\, T)$ s'annule transversalement.
Alors, nous noterons $A$ la sous-varit de $\R^m$ dfinie par $G=0$.

Il est bien connu que presque toute application polynomiale s'annule transversalement. On a donc le rsultat suivant.
\\
\\
\indent
Proposition.
Quel que soit le paramtre $T$, pour presque tout $G$ (parmi les applications polynomiales de degr $\leq D$), 
ce paramtre sera $G-$rgulier.

\subsubsection{La sous-varit $\overline{A}$}
Posons :
$$
 D_1(X,\, T) = \det \left( \frac{\partial G}{\partial X}(X,\, T) \circ \, \frac{{}^{t}\partial G}{\partial X}(X,\, T) \right) .
$$
Les points d'annulation $(X,\, T)$ de $D_1$ sont ceux o la diffrentielle partielle $\frac{\partial G}{\partial X}(X,\, T)$ n'est pas surjective. 

Posons :
\begin{eqnarray*}
D_2(X,\, T) &=& \det \left( \frac{\partial (f,\, G)}{\partial X}(X,\, T) \circ \, \frac{{}^{t}\partial (f,\, G)}{\partial X}(X,\, T) \right) 
\\ g (X,\, T) &=& D_2 (X,\, T) - \varepsilon^2 D_1 (X,\, T)
\end{eqnarray*}

Les deux types de points d'annulation $(X,\, T)$ de $g$ sont :

1)  Les points o la diffrentielle partielle $\frac{\partial G}{\partial X}(X,\, T)$ n'est pas surjective. 

2)  Parmi les points o cette diffrentielle partielle est surjective (et en lesquels la fibre de $G(.,\, T)$ est donc lisse),
ceux en lequels la norme du gradient de la restriction de $f$  la fibre de $G(.,\, T)$ vaut $\varepsilon$.
\\
\\
\indent
Posons :
$$ \overline{G} (X,\, T) = (g(X,\, T) ,\, G (X,\, T) ).$$

Un paramtre $T$ sera dit $ \overline{G}-$rgulier si l'application $ \overline{G}(.,\, T)$ s'annule transversalement.
Alors, nous noterons $\overline{A}$ la sous-varit de $\R^m$ dfinie par $ \overline{G}=0$.
\\
\\
\indent
Proposition.
Quel que soit le paramtre $T$, pour presque tous $G$, $f$ et $\varepsilon$, ce paramtre sera $ \overline{G}-$rgulier.
\\
\\
\indent
Dmonstration.
Soient $G$ s'annulant transversalement et $f$ quelconque. Il suffit de prendre $\varepsilon$
gal  une valeur rgulire $>0$ de la fonction qui associe  tout point de la sous-varit $A$
la norme du gradient en ce point de la restriction de $f$  $A$. D'aprs le thorme de Sard, presque tout $\varepsilon$ convient.

\subsubsection{La fonction de Morse $\|.\|^2$}
Pour tout $X\in \R^m$, on pose $\eta (X) = \|X\|^2$.
Un paramtre $G-$rgulier sera dit $(G,\,\eta)-$rgulier si la restriction de la fonction $\eta$
 la sous-varit $A$ est une fonction de Morse.
De mme, un paramtre $\overline{G}-$rgulier sera dit $(\overline{G},\,\eta)-$rgulier si la restriction de la fonction $\eta$
 la sous-varit $\overline{A}$ est une fonction de Morse.
\\
\\
\indent
Proposition.
Quel que soit le paramtre $T$, pour presque tout $G$, ce paramtre sera $(G,\, \eta)-$rgulier.
\\
\\
\indent
Dmonstration.
Le paramtre $T$ est suppos fix.
Notons ${\cal P}(G)$ la proprit  dmontrer.
Supposons $G$ choisie de telle sorte que $T$ soit $G-$rgulier.
Etant donn $C\in \R^m$, on posera $G_C(X,\, T) = G(X+C,\, T)$.

Il est bien connu qu'tant donne une sous-varit d'un espace euclidien, 
la fonction \og carr de la distance\fg\,  un point gnrique de l'espace est une fonction de Morse
sur la sous-varit.
Autrement dit, pour presque tout $C\in \R^m$, 
la proprit ${\cal P}(G_C)$ est vraie.
D'aprs le thorme de Fubini, ${\cal P}(G_C)$ est donc vraie pour presque tout couple $(G,\, C)$ (o $G$
dcrit l'espace des applications polynomiales de degr $\leq D$).

Remarquons que l'image rciproque d'un ensemble de mesure non nulle par la surjection $(G,\, C) \mapsto G_C$
est un ensemble de mesure non nulle. Donc la proprit ${\cal P}(G)$ est vraie pour presque tout $G$.
Le rsultat est dmontr.
\\
\\
\indent
Proposition.
Quel que soit le paramtre $T$, pour presque tous $G$, $f$ et $\varepsilon$, 
ce paramtre sera $(\overline{G},\, \eta)-$rgulier.
\\
\\
\indent
Dmonstration.
C'est presque pareil que la preuve prcdente.
On note ${\cal P}(G,\, f,\, \varepsilon)$ la proprit  dmontrer.
On pose $G_C(X,\, T) = G(X+C,\, T)$ et $f_C(X,\, T) = f(X+C,\, T)$.
Etc.

\subsubsection{Fubini}
Le thorme de Fubini permet bien sr de dduire des deux rsultats prcdents que pour presque tout $G$,
presque tout paramtre $T$ est $(G,\, \eta)-$rgulier
et que pour presque tous $G$, $f$ et $\varepsilon$,
presque tout paramtre $T$ est $(\overline{G},\, \eta)-$rgulier.

\subsubsection{Elimination}\label{elimination}

Il est vident que l'image d'un ensemble fini est un ensemble fini.
La version  paramtres est un peu moins vidente. En voici l'nonc :
\\
\\
\indent
Proposition.
Soit une application polynomiale $K: \R^m \times \R^p \rightarrow \R^m$.
Alors il existe une hypersurface $H \subset \R \times \R^p$ 
dfinie par un polynme non nul $h: \R \times \R^p \rightarrow \R$ de degr $\leq P(D)$
telle que, pour tout point $(X,\, T)$ en lequel
$K(X,\, T)$ est nul et  $\frac{\partial K}{\partial X}(X,\, T)$ est bijective, le point $(f(X,\, T),\, T)$
appartienne  $H$. 
\\
\\
\indent
Dmonstration.
Comme le nombre de composantes des applications polynomiales $K$, $f$ et $T$
est strictement suprieur au nombre des indtermines $X$, $T$, il existe une relation algbrique :
$$
h_1 ( K,\, f,\, T) = 0
$$ 
avec $h_1$ polynme non nul de degr $\leq P(D)$.
Les $m+1+p$ indtermines de $h_1$ seront notes $U$, $V$ et $T$.
On peut mettre $h_1$ sous la forme :
$$
h_1(U,\, V,\, T) = \sum_{\beta} U^{\beta} h_{\beta}(V,\, T).
$$
o $\beta \in \N^m$ dcrit un ensemble (non vide) d'exposants et o les $h_{\beta}(V,\, T)$ sont des polynmes non nuls.
Notons $\alpha$ le plus petit des $\beta$ pour l'ordre lexicographique et
vrifions que le polynme $h=h_{\alpha}$ convient.

Soit $(X_0,\, T_0)$ un point pour lequel $K(X_0,\, T_0)$ est nul et $\frac{\partial K}{\partial X}(X_0,\, T_0)$ est bijective.
Posons :
$$U=K(X,\, T).$$
Alors, d'aprs le thorme d'inversion locale, les $U$ et les $T$ forment un systme local de coordonnes
et les $X$ sont des fonctions analytiques des $U$ et des $T$.
$$
X = \varphi (U,\, T)
$$
Alors on peut crire :
$$
0 = h_1 ( U,\, f( \varphi (U,\, T),\, T),\, T) = \sum_{\beta} U^{\beta} h_{\beta} ( f( \varphi (U,\, T),\, T),\, T). 
$$ 
On simplifie par le facteur $U_1^{\beta_1}$ et on spcifie $U_1 = 0$
puis on simplifie par le facteur $U_2^{\beta_2}$ et on spcifie $U_2 = 0$ et ainsi de suite.
Cela dmontre la relation :
$$
0 = h_{\alpha} ( f( \varphi (0,\, T),\, T),\, T).
$$
Enfin on spcifie $T=T_0$ pour obtenir :
$$
0 = h_{\alpha} ( f( X_0,\, T_0),\, T_0). 
$$ 
Le rsultat est dmontr.

\subsubsection{Elimination et valeurs critiques d'une fonction de Morse}

Proposition.
Soit un polynme $\eta : \R^m \times \R^p \rightarrow \R$. 
Il existe une hypersurface $H \subset \R \times \R^p$ 
dfinie par un polynme non nul $h: \R \times \R^p \rightarrow \R$ de degr $\leq P(D)$
telle que, pour tout paramtre $G-$rgulier $T$ et tout point critique $X$ de Morse
de la restriction de $\eta(.,\, T)$  $A$,
le point $(f(X,\, T),\, T)$
appartienne  $H$. 
\\
\\
\indent
Dmonstration. 
Notons $I$ un $q-$uplet vrifiant $1 \leq i_1 < \dots < i_q \leq m$.
Travaillons sur l'ouvert $U_I$ o la diffrentielle partielle $\frac{\partial G}{\partial X_I}(X,\, T)$ est bijective.
Notons $J$ le complmentaire de $I$.

Afin d'appliquer la proposition (\ref{elimination}), choisissons l'application $K$ gale  :
$$\left( G,\, \frac{\partial \eta}{\partial X_J} - \frac{\partial \eta}{\partial X_I} \circ \left( \frac{\partial G}{\partial X_I} \right)^{-1}
\circ \frac{\partial G}{\partial X_J} \right) .$$
Alors on obtient un polynme $h_I$ qui convient sur l'ouvert $U_I$.
Comme il n'y a qu'un nombre fini de tels ouverts (et que ce nombre ne dpend que des dimensions), 
il suffit de prendre le produit $h$ des $h_I$ et le rsultat est dmontr.
\\
\\
\indent
Remarque. Le rsultat reste videmment vrai si on remplace $G$ par $\overline{G}$ et $A$ par $\overline{A}$.

\subsubsection{Marquage de certaines valeurs presque critiques}
Les donnes $G$, $f$, $\varepsilon$ et un paramtre $G-$rgulier $T$ tant fixes, 
on notera Crit l'ensemble des 
points $\varepsilon-$critiques de la restriction de $f$  $A$,
de norme $\leq 1$. Autrement dit :
$$ {\mbox Crit} = \{ X \in A\, ;\; g(X,\, T) \leq 0 \mbox{ et } \|X\| \leq 1\}.$$

Proposition.
Pour presque tous $f$, $G$ et $\varepsilon$, il existe une hypersurface $H \subset \R \times \R^p$ 
dfinie par un polynme non nul $h: \R \times \R^p \rightarrow \R$ de degr $\leq P(D)$
telle que, pour presque tout paramtre $G-$rgulier $T$, dans toute composante connexe de Crit il se trouve (au moins)
un point $X$ tel que l'image $f(X,\, T)$ appartienne  $H$.
\\
\\
\indent
Dmonstration. 
On choisit $f$, $G$ et $\varepsilon$ de telle sorte que presque tout paramtre $T$ soit
 la fois $(G,\, \eta)-$rgulier et $(\overline{G},\, \eta)-$rgulier.

Notons $C$ une composante connexe de Crit. Par compacit de $C$,
la restriction  $C$ de la fonction $\eta=\|.\|^2$ atteint son minimum 
en un point $X$. Distinguons deux cas.
\\
\\
\indent
1er cas : $X \in A \backslash \overline{A}$.
Alors le point $X$ est un minimum local de la restriction  $A$ de la fonction de Morse $\eta$.
Par limination, l'image $f(X,\, T)$ appartient  une hypersurface $H_1$.
\\
\\
\indent
2nd cas : $X \in \overline{A}$.
Alors le point $X$ est un minimum local de la restriction  $\overline{A}$ de la fonction de Morse $\eta$.
Par limination, l'image $f(X,\, T)$ appartient  une hypersurface $H_2$.
\\
\\
\indent
On conclut en prenant $H$ gale  la runion $H_1 \cup H_2$.

\subsubsection{Dmonstration du thorme (\ref{lezen 2}) dans le cas $n=1$}
Notons $H$ l'hypersurface donne par le rsultat prcdent qui assure que
pour presque tout paramtre $G-$rgulier $T$, dans toute composante $C$ de Crit,
il existe un $X\in C$ tel que $f(X,\, T)$ appartienne  $H$.
Il suffit donc de prouver que, pour tout $X' \in C$, on a $|f(X',\, T) - f(X,\, T)| < \varepsilon P(D)$.
C'est alors l'ingalit des accroissements finis qui permet de conclure. En effet,
d'une part, la norme du gradient de la restriction de $f$  $A$ est $\leq \varepsilon$ sur Crit
et d'autre part le diamtre par arcs d'une composante connexe $C$ de Crit est $\leq P(D)$
d'aprs le lemme 28 nonc dans \cite{Do99} (variante de la proposition 29 de \cite{Do96}).

Le cas $n=1$ de la rcurrence est donc dmontr, {\it a priori} pour presque tous $f$, $G$, $\varepsilon$ et $T$
mais, comme nous l'avons remarqu (cf. (\ref{gen gen})), a implique qu'il est vrai pour {\it tous}  $f$, $G$, $\varepsilon$ et $T$.

\subsection{Hrdit}

On suppose $n\geq 2$ et on note $f_1$, ..., $f_n$ les composantes de l'application $F$.

\subsubsection{Notations}
Etant donns un paramtre $T$ et un point $X$ d'annulation transversale de $G(.,\, T)$,
on notera $L(X,\, T)$ la diffrentielle en $X$ de la restriction de $F(.,\, T)$
au lieu d'annulation de $G(.,\, T)$.
On notera $l_i(X,\, T)$ la $i-$me composante de $L(X,\, T)$, c'est--dire la diffrentielle en $X$ de la restriction de $f_i(.,\, T)$
au lieu d'annulation de $G(.,\, T)$.  

\subsubsection{Gradient des composantes $f_i$}
Proposition.
Il existe une hypersurface $H_0 \in \R^n \times \R^p$, dfinie par un polynme non nul $h_0 :  \R^n \times \R^p \rightarrow \R$
de degr $\leq P(D)$ telle que pour tout paramtre $T\in \R^p$, si $X\in \R^m$ dsigne un point d'annulation transversale de $G(.,\, T)$ 
vrifiant $l_i(X,\, T) \leq \varepsilon$
pour (au moins) un $i$ compris entre $1$ et $n$,
alors il existe un point $Y\in \R^n$, $\varepsilon P(D)-$proche de $F(X,\, T)$ tel que $(Y,\, T)$ appartienne  $H_0$.
\\
\\
\indent
Dmonstration.
Pour chaque $i$, on applique le cas $n=1$ du thorme  $f_i$. On obtient donc une hypersurface de $\R \times \R^p$.
Son image rciproque dans $\R^n \times \R^p$ par la projection $(Y,\, T) \mapsto (Y_i,\, T)$ est une hypersurface $H_{0,i}$.
Alors la runion, quand $i$ varie entre $1$ et $n$, de ces hypersurfaces $H_{0,i}$ convient.

\subsubsection{Utilisation de l'hypothse de rcurrence}
On note $F_{\widehat{\i}} : \R^m \times \R^p \rightarrow \R^{n-1}$ l'application obtenue  partir de $F$
en gardant toutes les composantes sauf la $i-$me. 
En un point $X$ d'annulation transversale de $G(.,\, T)$, on note $L_{\widehat{\i}}(X,\, T)$
la restriction de $\frac{\partial F_{\widehat{\i}}}{\partial X} (X,\, T)$  $\ker \frac{\partial G}{\partial X} (X,\, T)$.
Un point $X$ d'annulation transversale de $G(.,\, T)$ sera dit $(G,\, i)-$rgulier si $l_i(X,\, T)$ n'est pas nul.
Si $X$ est $(G,\, i)-$rgulier, on note $\overline{L}_{\widehat{\i}}(X,\, T)$
la restriction de $L_{\widehat{\i}}(X,\, T)$  $\ker l_i(X,\, T)$.
\\
\\
\indent
Proposition. 
Il existe une hypersurface $H_i \in \R^n \times \R^p$, dfinie par un polynme non nul $h_i :  \R^n \times \R^p \rightarrow \R$
de degr $\leq P(D)$ telle que pour tout paramtre $T\in \R^p$, si $X\in \R^m$ dsigne un point $(G,i)-$rgulier 
en lequel module de surjectivit de $\overline{L}_{\widehat{\i}}(X,\, T)$ est  major par $\varepsilon$,
alors il existe un point $Y\in \R^n$, $\varepsilon P(D)-$proche de $F(X,\, T)$ tel que $(Y,\, T)$ appartienne  $H_i$.
\\
\\
\indent
Dmonstration.
Il suffit d'appliquer le thorme au rang $n-1$ aux donns suivantes :
\\
\\
\indent
1) On prend $m$ variables. Ce sont les $X$.

2) On prend $p+1$ paramtres. Ce sont les $T$ mais aussi une nouvelle indtermine que nous noterons $Y_i$.

3) On prend $q+1$ contraintes. Ce sont les composantes de $G(X,\, T)$ mais aussi $f_i(X,\, T) - Y_i$.

4) On prend comme application $F_{\widehat{\i}}$ dont l'espace d'arrive est bien de dimension $n-1$.
\\
\\
\indent
Alors on obtient, dans l'espace $\R^{n-1} \times \R^{p+1}$ identitfi  $\R^n \times \R^p$, une hypersurface $H_i$ qui convient. 

\subsubsection{L'hypersurface $H$}
Les deux rsultats prcdents fournissent, pour l'un, une hypersurface $H_0$
et, pour l'autre, $n$ hypersurfaces $H_i$. Notons $H$ la runion de ces $n+1$ hypersurfaces.  
Nous allons voir que cette hypersurface convient et le thorme sera donc dmontr. 

\subsubsection{Un peu d'algbre linaire}\label{MS decomp}

Rappelons que MS dsigne le {\it module de surjectivit.}
\\
\\
\indent
Lemme.
Soit une application linaire $u$ entre deux espaces euclidiens $E$ et $F$.
On suppose que $F$ est la somme directe orthogonale de deux sous-espaces $F_1$ et $F_2$.
Notons $u_1$ et $u_2$ les deux composantes de $u$ et notons $(u_2)_{\mbox{ker \it{u}}_1}$ la restriction de $u_2$
 $\mbox{ker } u_1$.

Alors l'ingalit suivante est vrifie :
$$
\mbox{MS } u_1 \mbox{ MS }(u_2)_{\mbox{ker \it{u}}_1} 
\leq \mbox{MS }u \; \left( \mbox{MS } u_1 + \mbox{MS } (u_2)_{\mbox{ker \it{u}}_1} + \| u_2 \|\right).
$$

Remarque. Ce lemme peut tre regard comme une version quantitative du fait que
 si $u_1$ et $(u_2)_{\mbox{ker \it{u}}_1}$ sont surjectives, $u$ le sera aussi.
\\
\\
\indent
Dmonstration du lemme.

Soit une forme linaire $\lambda = (\lambda_1,\, \lambda_2) \in F^* = F_1^* \oplus F_2^*$. La relation :
$$ 
\lambda \circ u = \lambda_1 \circ u_1 + \lambda_2 \circ u_2
$$
implique, par ingalit triangulaire :
\begin{eqnarray*}
\| \lambda_1 \| \mbox{ MS }u_1 & \leq & \| \lambda_1 \circ u_1 \|
\\ & \leq & \| \lambda \circ u \| + \| \lambda_2 \| \| u_2 \|.
\end{eqnarray*}
Par ailleurs, on a :
\begin{eqnarray*}
\| \lambda_2 \| \mbox{ MS }(u_2)_{\mbox{ker \it{u}}_1} & \leq & \left\| (\lambda_2 \circ u_2)_{\mbox{ker \it{u}}_1} \right\|
\\ & = &  \left\| (\lambda \circ u)_{\mbox{ker \it{u}}_1} \right\|
\\ & \leq & \| \lambda \circ u \|.
\end{eqnarray*}
Appliquons  nouveau l'ingalit triangulaire :
$$
\| \lambda \| \leq \| \lambda_1 \| + \| \lambda_2\|.
$$
Rassemblons les trois ingalits prcedentes :
\begin{eqnarray*}
\| \lambda \| \; \mbox{MS } u_1 \mbox{ MS }(u_2)_{\mbox{ker \it{u}}_1} 
& \leq & ( \| \lambda_1 \| + \| \lambda_2\|) \; \mbox{MS } u_1 \mbox{ MS }(u_2)_{\mbox{ker \it{u}}_1} 
\\ & \leq & ( \| \lambda \circ u\| + \| \lambda_2 \| \| u_2\| + \| \lambda_2\| \, \mbox{MS } u_1 )\, \mbox{ MS }(u_2)_{\mbox{ker \it{u}}_1} 
\\ & \leq & \| \lambda \circ u\| \; ( \mbox{ MS }(u_2)_{\mbox{ker \it{u}}_1}  +\| u_2\| + \mbox{MS } u_1 )
\end{eqnarray*} 
On conclut en appliquant cette ingalit  une forme $\lambda$ de norme $1$ ralisant le minimum de $\| \lambda \circ u \|$ 
(c'est--dire vrifiant : $\mbox{MS }u = \| \lambda \circ u\|$ ).
\\
\\
\indent
Proposition.
Soit $X$ un point $(G,\, i)-$rgulier. 
Alors, l'ingalit suivante est vrifie :
$$
\| l_i(X,\, T) \| \mbox{ MS }(\overline{L}_{\widehat{\i}}(X,\, T))
\leq \mbox{MS }(L(X,\, T)) \; \left( \| l_i(X,\, T) \| + \mbox{MS } (\overline{L}_{\widehat{\i}}(X,\, T)) + \| L_{\widehat{\i}}(X,\, T)) \|\right).
$$

Dmonstration. 
Il suffit d'appliquer le lemme en remarquant que $ l_i(X,\, T) $ prend ses valeurs dans $\R$
et que son module de surjectivit est donc simplement gal  sa norme $\|  l_i(X,\, T) \| $. 

\subsubsection{Fin de la dmonstration des thormes (\ref{lezen}) et (\ref{lezen 2})}
Soit $X$ un point d'annulation transversale de $G(.,\, T)$.
Choisisssons un $i$ compris entre $1$ et $n$ qui maximise la norme de $l_i(X,\, T)$.
Alors l'identit de Pythagore implique :
$$
\| L_{\widehat{\i}}(X,\, T) \|  \leq \sqrt{n-1} \;  \| l_i(X,\, T) \|. 
$$

Notons $H(T)$ (resp. $H_0(T)$ et $H_i(T)$) l'ensemble des $Y\in \R^n$ pour lesquels $(Y,\, T)$ appartient  $H$ (resp.  $H_0$ et  $H_i$).
On suppose que $F(X,\, T)$ n'est \og pas trop proche \fg\, de $H(T)$
(en un sens que nous allons prciser). Alors $F(X,\, T)$ n'est \og pas trop proche \fg\, de $H_0(T)$ et l'ingalit suivante est donc vrifie :
$$\| l_i(X,\, T) \|  \geq \varepsilon.$$
Notamment, le point $X$ est $(G,\, i)-$rgulier. 
De mme, le point $F(X,\, T)$ n'est \og pas trop proche \fg\, de $H_i(T)$ et l'ingalit suivante est donc vrifie :
$$\mbox{MS} (\overline{L}_{\widehat{\i}}(X,\, T)) \geq \varepsilon.$$
Afin de s'assurer que toutes ces ingalit soient vraies, il suffit de choisir un polynme $P$ convenable 
et de dire que l'expression  \og $F(X,\, T)$ n'est pas trop proche de $H(T)$\fg \, signifie 
qu'il n'est $\varepsilon P(D)-$proche d'aucun point de $H(T)$.

La proposition (\ref{MS decomp}) permet d'crire :
$$
\| l_i(X,\, T) \| \mbox{ MS }(\overline{L}_{\widehat{\i}}(X,\, T))
\leq \mbox{MS }(L(X,\, T)) \; \left( (1+\sqrt{n-1}) \| l_i(X,\, T) \| + \mbox{MS } (\overline{L}_{\widehat{\i}}(X,\, T)) \|\right).
$$

On en dduit :
\begin{eqnarray*}
\mbox{MS} (L(X,\, T)) &\geq &  (2+\sqrt{n-1})^{-1} \times \min \left\{ \| l_i(X,\, T) \| ,\, \mbox{MS } (\overline{L}_{\widehat{\i}}(X,\, T)) \right\}
\\ & \geq &  (2+\sqrt{n-1})^{-1} \times  \varepsilon.
\end{eqnarray*}
Bien sr, on peut faire disparatre la constante $(2+\sqrt{n-1})^{-1}$ quitte  changer de polynme $P(D)$.

Ouf ! Le thorme (\ref{lezen 2}) est dmontr ainsi que le thorme (\ref{lezen}) qui en est un cas particulier.

\section{Transversalit quantitative pour les applications polynomiales et holomorphes}

\subsection{Volume du voisinage d'une hypersurface}\label{vol hyp}
Proposition.
Soit une hypersurface $H$ dfinie dans $\R^n$ par une fonction polynomiale non nulle $f:\R^n \rightarrow \R$.
Notons $B$ la boule unit de $\R^n$ et posons : 
$$V_{\varepsilon}=\left\{ x \in \R ^n ;\; \exists y \in H \;\; \|x-y\|\leq \varepsilon \right\}$$
pour tout $\varepsilon >0$.

Alors le volume de $B\cap V_{\varepsilon}$ sera major par :
$$
\varepsilon \times P(\mbox{deg }f)
$$
o $P(\mbox{deg }f)$ dsigne une fonction polynomiale du degr de $f$  (la dimension $n$ tant fixe).
\\
\\
\indent
Dmonstration.
L'ensemble $B\cap V_{\varepsilon}$ est form de deux types de points :

1. Certains sont $\varepsilon-$proches du bord de $B$ et leur volume est major par le produit d'$\varepsilon$ et d'une constante (le volume $(n-1)-$dimensionnel de la sphre unit).

2. Les autres sont $\varepsilon-$proches de $B\cap  H$. Ils sont donc $2\varepsilon-$proches d'une $\varepsilon-$discrsitsation $F$.
Le volume tudi est donc major par le produit du cardinal de $F$ et du volume d'une boule de rayon $2\varepsilon$.
Le volume d'une telle boule est bien sr de l'ordre de $\varepsilon^n$ et,
d'aprs la proposition (\ref{comp hyp}), le cardinal d'une discrtisation admet une majoration du type : $\varepsilon^{1-n} \times P($deg $f)$.
Cela donne la majoration recherche.

\subsection{Point pas trop proche d'un hypersurface}\label{pas proche}
Proposition.
Avec les mmes notations, il existe un point de $\R^n$ de norme $\leq \varepsilon$ dont la distance  $H$ soit suprieure  :
$$\frac{\varepsilon}{P(\mbox{deg }f)}$$
o $P$ dsigne une fonction polynomiale  valeurs strictement positives (la dimension $n$ tant fixe).
\\
\\
\indent
Dmonstration.
On se ramne par une homothtie au cas $\varepsilon =1$.
 Ensuite d'aprs la proposition (\ref{vol hyp}),
le volume de l'ensemble $V_{\eta}$ des  points $\eta -$proches de $H$
est major par $\eta \times P($deg $f)$ pour un certain polynme $P$.
Quitte  changer de polynme $P$, si on prend $\eta$
gal  $(P($deg $f))^{-1}$ le volume de $V_{\eta}$ sera strictement infrieur  une constante donne,
par exemple gale au volume de la boule unit. Alors, 
par un argument vident de volume, la boule unit n'est pas incluse dans $V_{\eta}$.

\subsection{Transversalisation quantitative d'une application polynomiale}\label{trans quant pol}
Proposition. 

Soit une application polynomiale $f$ entre deux espaces euclidiens $X$ et $Y$.
Pour tout rel $\varepsilon >0$, il existe une \og bonne\fg\; valeur rgulire, c'est--dire
un point $y\in Y$ de norme $\leq \varepsilon$ tel que, pour tout point $x\in X$ de norme $\leq 1$,
le module de transversalit en $x$ de l'application $f - y$ vrifie la minoration suivante :
$$
\mbox{MT}(f-y, \, x) \geq \frac{\varepsilon }{P(\mbox{deg }f)}
$$ 
o $P(\mbox{deg }f)$ dsigne une fonction polynomiale de deg $f$ (les dimensions de $X$ et $Y$ tant fixes).
\\
\\
\indent
Dmonstration.
Notons $H$ l'hypersurface fournie par le thorme (\ref{lezen}). 
En appliquant la proposition (\ref{pas proche})  $H$, on obtient un point $y$ 
qui convient car il est $\varepsilon' -$distant de $f(x)$, pour tout point $\varepsilon' -$critique $x$ de norme $\leq 1$,
le rel $\varepsilon'$ pouvant tre choisi de telle sorte que le
quotient $\frac{\varepsilon}{\varepsilon'}$ soit polynomial en le degr de $H$, lui-mme polynomial en le degr de $f$. 

\subsection{Un rsultat lmentaire}\label{lem elem}
Lemme.
Soient deux rels $C$ et $D>0$ et une fonction polynomiale $M_1: \R \rightarrow \R_+^*$.
Alors il existe un rel $\varepsilon_0 >0$ et une fonction polynomiale $M_2 : \R \rightarrow \R_+^*$ tels que, pour tout rel $\varepsilon \in ]0,\, \varepsilon_0 [$, il existe un entier $n\geq 0$ vrifiant :
$$
\frac{\varepsilon }{M_1(n)} -Ce^{-Dn} \geq \frac{\varepsilon}{M_2\left(\log \frac{1}{\varepsilon}\right)}\;\;.
$$
\\
\indent
Dmonstration.
Notons $n$ le plus petit entier $\geq \frac{2}{D} \log \frac{1}{\varepsilon}$. Alors : 

$$Ce^{-Dn} \leq C\varepsilon ^2  \leq \frac{\varepsilon}{2M_1\left(\log \frac{1}{\varepsilon}\right)}$$ pour $\varepsilon$ assez petit.
Par ailleurs, pour cette valeur de $n$ :
$$2M_1(n) \leq M_2\left(\log \frac{1}{\varepsilon}\right)$$ avec un $M_2$ polynomial.

On en dduit :
\begin{eqnarray*}
\frac{\varepsilon }{M_2\left(\log \frac{1}{\varepsilon}\right)} & \leq &
\frac{\varepsilon}{2M_1(n)} \\ & \leq &
\frac{\varepsilon }{M_1(n))} -Ce^{-Dn}\;\;.
\end{eqnarray*}

\subsection{Transversalisation quantitative d'une application holomorphe}\label{trans hol}
Thorme.

Soient deux espaces vectoriels hermitiens $X$ et $Y$, un sous $\R-$espace vectoriel $A \subset X$
et deux rels $0<r<R$. Pour tout rel $s>0$, on notera $B(s)$ la boule centre en l'origine de $X$ de rayon $s$.
Soit une application holomorphe $f: B(R) \rightarrow Y$ vrifiant $\| f(x) \| \leq 1$ en tout point $x\in B(R)$.
Alors, pour tout $\varepsilon >0$ suffisament petit, il existe une  \og bonne\fg\; valeur rgulire, c'est--dire
un point $y\in Y$ de norme $\leq \varepsilon$ tel que, pour tout point $x\in A\cap B(r)$,
le module de transversalit en $x$ de la restriction $(f-y)_A$ au sous-espace $A$ de l'application $f - y$ vrifie la minoration suivante :
$$
\mbox{MT}((f-y)_A, \, x) \geq \frac{\varepsilon}{P\left( \log \frac{1}{\varepsilon}\right)}
$$ 
o $P(\log \frac{1}{\varepsilon})$ dsigne une fonction polynomiale de $\log \frac{1}{\varepsilon}$ (les dimensions et les rels $R$ et $r$ tant fixs)
et o l'expression \og pour tout $\varepsilon$ suffisament petit \fg\; signifie pour tout $\varepsilon$ major par un certain rel $>0$
ne dpendant que des dimensions et des rels $R$ et $r$.
\\
\\
\indent
Dmonstration.
Les polynmes de Taylor $(p_n)$ de l'application holomorphe $f$ convergent vers $f$ sur la boule $B(r)$  une vitesse exponentielle (voir \cite{HeLe84} p.56). 
C'est notamment vrai en norme ${\cal C}^1$ :
$$
\| f-p_n \|_{{\cal C}^1(B(r))} \leq Ce^{-Dn}
$$
pour des rels $C$ et $D>0$. Par ailleurs, la proposition (\ref{trans quant pol}) fournit un point $y\in Y$ norme $\leq \varepsilon$, tel que 
le module de transversalit de la restriction de $p_n-y$  $A\cap B(R)$ admette une minoration du type $\varepsilon \times (P(n))^{-1}$.
Remarquons que le module de transversalit est videmment $1-$lipschitzien pour la norme ${\cal C}^1$, on en dduit que celui de la restriction de $f-y$
est minor par : $\varepsilon \times (P(n))^{-1} -Ce^{-Dn}$.
Alors, le lemme (\ref{lem elem}) permet choisir une valeur convenable de $n$ qui permette de conclure.
\\
\\
\indent
Nous aurons besoin du rsultat suivant, qui se dduit du thorme prcdent par un simple changement d'chelle :
\\
\\
\indent
Corollaire.

On suppose donns un entier $k\geq 1$ et
une application holomorphe $f: B(k^{-\frac{1}{2}}R) \rightarrow Y$ vrifiant $\| f(x) \| \leq 1$ en tout point $x\in B(k^{-\frac{1}{2}}R)$.
Alors, pour tout $\varepsilon >0$ suffisament petit, il existe une  \og bonne\fg\; valeur rgulire, c'est--dire
un point $y\in Y$ de norme $\leq \varepsilon$ tel que, pour tout point $x\in A\cap B(k^{-\frac{1}{2}}r)$,
le module de transversalit en $x$ de poids $(1,\, k^{-\frac{1}{2}})$ de la restriction $(f-y)_A$ au sous-espace $A$ 
de l'application $f - y$ vrifie la minoration suivante :
$$
\mbox{MT}\left( (f-y)_A, \, x;\; 1,\, k^{-\frac{1}{2}}\right) \geq \frac{\varepsilon}{P\left( \log \frac{1}{\varepsilon}\right)}
$$ 
o la fonction polynomiale $P$ et le majorant de $\varepsilon$ (implicite dans l'expression \og $\varepsilon $ suffisament petit\fg )
sont ceux du thorme.
Insistons notamment sur le fait qu'ils ne dpendent pas de l'entier $k$.

\section{Construction de beaucoup de sections d'un fibr trs positif}
\subsection{Section concentre prs d'un point}\label{section concentree}

Un fibr trs positif admet des sections approximativement holomorphes concentres.
Cette ide se traduit par des estimes.

Pour tout point $x_0$ de la varit presque khlrienne $X$, on notera $L_{x_0}^{-1}\otimes L$ le fibr en droites hermitiennes dont
la fibre au-dessus d'un point $x\in X$ est $L_{x_0}^{-1}\otimes L_x$.

Petite subtilit : la varit $X$ n'est pas suppose compacte. Afin de pouvoir crire des estimes uniformes en $x_0$, il nous
faudra astreindre $x_0$  rester dans une partie compacte. Pour cette raison, nous prendrons en fait $x_0$
dans la sous-varit $Y$ qui, elle, est suppose compacte.

Le rsultat principal de ce paragraphe est qu'il existe une section $c_{x_0}$ de $L_{x_0}^{-1}\otimes L$ que l'on peut choisir de
telle sorte que soient vrifies les quatre estimes que nous allons bientt noncer. 
Nous ne donnerons pas d'nonc plus formel :
nous renvoyons  \cite{Do96} et \cite{Au97} pour la construction de cette section ainsi que pour la dmonstration des quatre estimes
que nous nous contenterons d'noncer.

Les deux premires estimes portent sur les sections \og concentres\fg, c'est--dire sur les puissances positives $c_{x_0}^k$ avec $k\geq 1$.
La puissance $c_{x_0}^k$ est bien sr une section du fibr $L_{x_0}^{-k}\otimes L^k$. Sur ce fibr, nous utiliserons la mtrique et la connexion induites par celles de $L$.

La premire estime est :
\begin{eqnarray}\label{premiere estimee}
\left\| (\nabla^m c_{x_0}^k)_x \right\| 
& \leq  &
k^{\frac{m}{2}} \times P\left(k^{\frac{1}{2}}d(x_0,\, x)\right) \times \mbox{exp}\left(-\frac{k\pi}{2} d^2(x_0,\, x)\right)
\end{eqnarray}
o $x_0$ dsigne un point de $Y$, o $x$ dsigne un point de $X$, o l'ordre de drivation $m$ est un entier $\geq 0$ quelconque,
o l'exposant $k$ est un entier $\geq 1$ quelconque et o $P$ dsigne une fonction polynomiale qui dpend 
de la varit presque khlrienne $X$, de sa pr-quantification $L$, de la sous-varit $Y$ et de l'ordre de drivation $m$.
Remarquons que $P$ ne dpend ni de l'exposant $k$ ni des points $x_0$ et $x$.

La seconde estime est : 
\begin{eqnarray}
\left\| (\nabla^m \overline{\partial} c_{x_0}^k)_x \right\| 
& \leq  &
k^{\frac{m}{2}} \times P\left(k^{\frac{1}{2}}d(x_0,\, x)\right) \times \mbox{exp}\left(-\frac{k\pi}{2} d^2(x_0,\, x)\right)
\end{eqnarray}
avec les mmes conventions.

La troisime et la quatrime estime porteront sur les puissances ngatives $c_{x_0}^{-k}$ avec $k\geq 1$.
Bien sr, ces puissances sont dfinies l o $c_{x_0}$ ne s'annule pas.

La section $c_{x_0}$ ne s'annule pas sur la boule $B(x_0,\, \alpha)$ pour un certain rel $\alpha >0$ uniforme en $x_0$. En particulier,
tant donn un rel $R>0$, la section $c_{x_0}$ ne s'annulera pas sur la boule $B(x_0,\, R\times k^{-\frac{1}{2}})$ si l'exposant $k$ 
est assez grand. Ici, \og assez grand \fg \, signifie que $k$ est suprieur  une borne dpendant de $X$, $L$, $Y$ et $R$.

La troisime et la quatrime estime seront valables pour $k$ assez grand (c'est--dire suprieur  une borne dpendant de $X$, $L$, $Y$ et $R$)
et pour $x$ assez proche de $x_0$ (c'est--dire $x\in B(x_0,\, R\times k^{-\frac{1}{2}})$). Elles s'crivent :

\begin{eqnarray}
\|(\nabla^m c_{x_0}^{-k})_x\| & \leq & C\, k^{\frac{m}{2}}
\\
\|(\nabla^m \overline{\partial} c_{x_0}^{-k})_x\| & \leq & C\, k^{\frac{m}{2}}
\end{eqnarray}
o $C$ dsigne un rel qui dpend de $X$, $L$, $Y$, $m$ et $R$.

\subsection{Taille d'une discrtisation.}\label{taille discr}
Proposition.
Soient une varit riemannienne $X$ et une partie compacte $Y$ de $X$.
Il existe un rel $\varepsilon_0>0$ et une fonction polynomiale $P$ vrifiant la proprit suivante :
\\
\\
\indent
Soient deux rels $\varepsilon \in ]0,\, \varepsilon_0 [$ et $C>0$ et une partie finie $F$ de $Y$ satisfaisant, pour tous $x$, $y\in F$, la condition d'espacement : $d(x,y)\geq \varepsilon$.
Alors, pour tout $x\in X$, le cardinal de l'intersection $F \cap B(x,\, C)$ est major par $P\left( \frac{C}{\varepsilon} \right)$. 
\\
\\
\indent
Dmonstration. 
Comme $Y$ est compacte, il existe des rels strictement positifs $\alpha$ et $R_1$ tels que tout point $x_0 \in Y$ et tout rel $R\in]0, \, R_1]$ vrifient :
$$
\mbox{Vol }B(x_0,\, R) \geq \alpha R^d
$$
o $d$ dsigne la dimension de $X$.

Nous appellerons rayon d'injectivit relatif, que nous noterons $R(X,\, Y)$, le minimum des rayons d'injectivit de $X$ en les points de $Y$.
Il est strictement positif par compacit de $Y$. Posons $R_2 = \frac{1}{2}R(X,\, Y)$ et :
$$
V(Y, \, R_2)= \{ x\in X ;\;\; \exists y\in Y \;\; d(x,y)\leq R_2  \}.
$$
L'ingalit $R_2 < R(X,\, Y)$ implique que $V(Y, \, R_2)$ est compact. 
Il existe donc un rel strictement positif $\beta$ tel que tout point $x\in X$ et tout rel $R> 0$ vrifient : 
$$
\mbox{Vol }(B(x, \, R) \cap V(Y, \, R_2)) \leq \beta R^d.
$$
Posons $\varepsilon_0 = 2 \min \{R_1,\, R_2\}$. Alors, tout $\varepsilon \leq \varepsilon_0$ satisfait les deux conditions suivantes, valables pour tout $x_0 \in Y$ :
\begin{eqnarray*}
\mbox{Vol }B\left( x_0,\, \frac{\varepsilon}{2} \right) & \geq & \alpha \left( \frac{\varepsilon}{2}\right)^d
\\
B\left( x_0,\, \frac{\varepsilon}{2} \right) & \subset & V(Y, \, R_2).
\end{eqnarray*} 
Si le point $x_0$ appartient de plus  la boule $B(x, \, C)$ (pour un $x \in X$), l'ingalit triangulaire impliquera :
\begin{eqnarray*}
B\left( x_0,\, \frac{\varepsilon}{2} \right) & \subset &B\left( x, \, C+\frac{\varepsilon}{2} \right) \cap V(Y, \, R_2).
\end{eqnarray*} 

Enfin la condition d'espacement implique que les boules de rayon $\frac{\varepsilon}{2}$ centres en deux points de $F$ distincts seront disjointes, ce qui permet de conclure :
\begin{eqnarray*}
\sum_{x_0} \mbox{Vol }B\left( x_0,\, \frac{\varepsilon}{2} \right) & \leq & \mbox{Vol }\left( B\left( x,\, C+\frac{\varepsilon}{2} \right)\cap V(Y, \, R_2)\right)
\end{eqnarray*}
o $x_0$ dcrit $F\cap B(x,\, C)$.
Si on note $n$ le cardinal de $F\cap B(x,\, C)$, on aura donc :
\begin{eqnarray*}
n \, \alpha \left( \frac{\varepsilon}{2}\right)^d & \leq & \beta \left( C+\frac{\varepsilon}{2} \right)^d
\end{eqnarray*}
Autrement dit :
\begin{eqnarray*}
n & \leq & \frac{\beta}{\alpha} \left( \frac{2C}{\varepsilon} +1 \right)^d.
\end{eqnarray*}
Le membre de droite de cette dernire ingalit est bien polynomial en $\frac{C}{\varepsilon} $.

\subsection{Un rsultat lmentaire.}\label{element}
Lemme.
Soient un rel $C>0$ et un polynme $M_1 : \R \rightarrow \R$. Alors il existe un polynme $M_2 : \R \rightarrow \R$ vrifiant, pour tout entier $N\geq 0$ :
$$
\sum_{n \geq N} M_1(n) \times \mbox{exp}(-Cn) = M_2(N) \times \mbox{exp}(-CN).
$$
\\
\\
\indent
Dmonstration. 
En associant  une suite $u=(u_n)_n$ la suite $\phi(u)$ dfinie par :
$$
\phi (u)_n = u_n -u_{n+1},
$$
on dfinit un endomorphisme $\phi$ de l'espace des suites de la forme suivante : $(M(n) \times \mbox{exp}(-Cn))_n$ avec $M$ polynme de degr $\leq$ deg $M_1$. 
Comme $\phi$ est injectif (et que la dimension de l'espace est finie), $\phi$ sera surjectif.

\subsection{Combinaisons linaires de sections concentres.}\label{combi}
Retournons au contexte qui nous occupe : on suppose donnes une varit presque khlrienne $X$
munie d'une pr-quantification $L$ et une sous-varit compacte $Y$ ainsi qu'un rang entier $r\geq 1$. Les sections $c_{x_0}$ sont 
les mmes que prcdemment.
\\
\\
\indent
Proposition.
Soit une partie finie $F$ de $Y$ satisfaisant, pour tous $x$, $y\in F$, la condition d'espacement :
$d(x,\, y) \geq k^{-\frac{1}{2}}$ pour un certain entier $k\geq 1$. Supposons donn, pour tout point $x_0 \in F$,
un lment $\alpha^{x_0} \in rL_{x_0}^k$ de norme
$\leq 1$.
Notons $s$ la section du fibr vectoriel $rL^k$ dfinie par  :
$$
s=\sum_{x_0 \in F} \alpha^{x_0} c_{x_0}^k.
$$
Pour tout $x\in X$, posons :
$$
d(x,\, F) = \min_{y\in F} d(x,\, y).
$$
Alors les deux majorations suivantes sont vrifies :
\begin{eqnarray*}
\left\| (\nabla^m s)_x \right\| 
& \leq  &
C \times k^{\frac{m}{2}} \times \mbox{exp}\left(-\frac{k\pi}{3} d^2(x,\, F)\right)
\\
\left\| (\nabla^m \overline{\partial} s)_x \right\| 
& \leq  &
C \times k^{\frac{m}{2}} \times \mbox{exp}\left(-\frac{k\pi}{3} d^2(x,\, F)\right)
\end{eqnarray*}
o $x$ dsigne un point de $X$, o $m$ dsigne un entier $\geq 0$
et o $C$ dsigne un rel qui dpend de $X$, $L$, $Y$ et $m$.
\\
\\
\indent
Dmonstration.
Une partition de $F$ est donne par les sous-ensembles suivants :
$$
F_n = \left\{ y \in F ; \;\; n \leq k d^2(x,\, F) <n+1 \right\}
$$
pour $n\geq$ Ent$(k d^2(x,\, F))$, o Ent dsigne la partie entire.
\begin{eqnarray*}
\mbox{card }F_n & \leq & \mbox{card }\left\{ y \in F ; \;\; k d^2(x,\, F) <n+1 \right\}
\\ & \leq & P_1 \left( \frac{k^{-\frac{1}{2}}\sqrt{n+1}}{k^{-\frac{1}{2}}}\right) = P_1 (\sqrt{n+1})
\end{eqnarray*}
o $P_1 $ dsigne le polynme donn par la proposition (\ref{taille discr}).

Soient deux point $x_0 \in F$ et $x\in X$ qui vrifient : $\sqrt{n}\leq k^{\frac{1}{2}} d(x_0,\, x) < \sqrt{n+1}$. 
Renommons $P_2$ le polynme $P$ intervenant dans l'estime (\ref{premiere estimee}) du (\ref{section concentree}). Sans perte de gnralit, on peut supposer $P_2$ croissant. 
Alors : $$P_2(k^{\frac{1}{2}} d(x_0,\, x)) \leq P_2(\sqrt{n+1})$$
et par ailleurs :
$$\mbox{exp}\left(-\frac{k\pi}{2} d^2(x_0,\, x)\right) \leq \mbox{exp}\left(-\frac{\pi n}{2}\right)$$
et donc, d'aprs l'estime (\ref{premiere estimee}) du (\ref{section concentree}) :
\begin{eqnarray*}
\left\| \alpha^{x_0}(\nabla^m c_{x_0}^k)_x \right\| 
& \leq & 
k^{\frac{m}{2}} \times \| \alpha^{x_0} \| \times P_2(\sqrt{n+1}) \times \mbox{exp}\left(-\frac{\pi n}{2}\right)
\\ & \leq & 
k^{\frac{m}{2}} \times P_2(\sqrt{n+1}) \times \mbox{exp}\left(-\frac{\pi n}{2}\right).
\end{eqnarray*}
Donc :
\begin{eqnarray*}
\sum_{x_0 \in F_n} \left\| \alpha^{x_0}(\nabla^m c_{x_0}^k)_x  \right\| 
& \leq & \mbox{card }F_n  \times k^{\frac{m}{2}} \times P_2(\sqrt{n+1}) \times \mbox{exp}\left(-\frac{\pi n}{2}\right)
\\ & \leq & k^{\frac{m}{2}} \times P_1(\sqrt{n+1}) \times P_2(\sqrt{n+1}) \times \mbox{exp}\left(-\frac{\pi n}{2}\right)
\\ & \leq & k^{\frac{m}{2}} \times P_3(n) \times \mbox{exp}\left(-\frac{\pi n}{2}\right)
\end{eqnarray*}
o $P_3$ dsigne un polynme.
\begin{eqnarray*}
\left\| (\nabla^m s)_x \right\| 
 & \leq & k^{\frac{m}{2}} \times \sum_{n\geq \mbox{ Ent} (k d^2(x,\, F))} P_3(n) \times \mbox{exp}\left(-\frac{\pi n}{2}\right)
\\  & \leq & k^{\frac{m}{2}} \times P_4\left( \mbox{Ent} (k d^2(x,\, F))\right) \times \mbox{exp}\left(-\frac{\pi }{2}\mbox{Ent} (k d^2(x,\, F))\right)
\end{eqnarray*}
o $P_4$ dsigne le polynme fourni par le lemme (\ref{element}).
Les croissances compares de l'exponentielle et du polynme permettent d'en
dduire une majoration de $\left\| (\nabla^m s)_x \right\| $ du type recherch. La majoration de $\left\| (\nabla^m \overline{\partial} s)_x \right\| $ s'obtient de faon analogue.

\section{Procd de transversalisation de Donaldson}

Ce dernier chapitre achve la dmonstration du thorme de Donalson-Auroux relatif
en exposant le procd de transversalisation de Donaldson. 
L encore, le cas relatif est trs semblable au cas absolu. 

\subsection{Lemme $\overline{\partial}$ et changement d'chelle}\label{d bar ech}
Notons $E$ et $F$ deux espaces vectoriels hermitiens.
Pour tout rel $s>0$, on notera $B(s)$ la boule centre en l'origine de $E$ de rayon $s$. 
\\
\\
\indent
Proposition.
Soit une application holomorphe $f: B(k^{-\frac{1}{2}}R) \rightarrow F$, pour un certain entier $k\geq 1$ et un certain rel $R>0$. 
On suppose que $f$ vrifie la majoration suivante : 
\begin{eqnarray*}
\sup_{x\in B(k^{-\frac{1}{2}}R)} \| D^m\overline{\partial}f(x) \| & \leq & k^{\frac{m}{2}},
\end{eqnarray*}
pour tout entier $m$ compris entre $0$ et un certain entier $m_{max}$.
Soit un rel $R_1$ vrifiant $0<R_1<R$.

Alors il existe une application holomorphe $g: B(k^{-\frac{1}{2}}R_1) \rightarrow F$ vrifiant la majoration suivante : 
\begin{eqnarray*}
\sup_{x\in B(k^{-\frac{1}{2}}R_1)} \| D^mf(x)-D^mg(x) \| & \leq & C\, k^{\frac{m-1}{2}}
\end{eqnarray*}
pour tout $m$ compris entre $0$ et $m_{max}$, le rel $C$ dpendant des dimensions de $E$ et de $F$, des rels $R_1$ et $R$ et de l'entier $m_{max}$.
\\
\\
\indent
Dmonstration. Le cas particulier $k=1$ est le lemme $\overline{\partial}$ usuel sur une boule et le cas gnral s'y ramne en composant  la source avec une homothtie.
\\
\\
\indent
Remarque. Il est bien connu que le lemme $\overline{\partial}$ sur la boule est, en fait, un thorme difficile. Pour la dmonstration de ce \og lemme \fg, on renvoie  \cite{HeLe84} (p. 59-60).

\subsection{Transversalisation prs d'un point}\label{trans point}
Proposition.
Soit $s$ une section $(1,\, 2)-$contrle et $(1,\, 1)-$approximativement holomorphe  du fibr vectoriel $rL^k$, pour un certain entier $k$ assez grand.
Notons $R$ un rel $>0$ quelconque, $\varepsilon >0$ un rel assez petit et $y_0$ un point de la sous-varit $Y$.

Alors il existe un lment $\alpha \in rL_{y_0}^k$ de norme $\leq \varepsilon$ tel que si
on pose $s'=s+\alpha \, c_{y_0}^k$ et qu'on note $s_Y'$ la restriction de $s'$  la sous-varit $Y$, 
le module de transversalit pondre $\mbox{MT}(s_Y', \, y;\; 1,\, k^{-\frac{1}{2}})$ sera minor par $\varepsilon \times P \left( \log \frac{1}{\varepsilon}\right)^{-1}$ pour tout $y \in Y \cap B(y_0, \, R\times k^{-\frac{1}{2}})$.

Prcisons que \og $\varepsilon$ assez petit \fg\, signifie que $\varepsilon$ est major par un rel $>0$ 
qui dpend des donnes $X$, $L$, $Y$, $r$ et $R$
et que $P$ dsigne une fonction polynomiale  valeurs strictement positives qui dpend aussi de ces donnes.
Enfin \og $k$ assez grand \fg\, signifie que $k$ est suprieur  une borne qui dpend de $\varepsilon$ et
des donnes $X$, $L$, $Y$, $r$ et $R$. 
\\
\\
\indent
Dmonstration.
Avant de commencer cette dmonstration, il faut avertir le lecteur d'un choix didactique d'exposition qui consiste  ne pas supposer la section $(1,\, 2)-$contrle et $(1,\, 1)-$approximativement holomorphe mais $(1,\, m_1)-$contrle et $(1,\, m_2)-$approximativement holomorphe pour des entiers $m_1 \geq 2$ et $m_2\ge 1$. Bien sr, a ne change rien au contenu mathmatique de l'nonc mais a permet peut-tre une prsentation plus uniforme des calculs.

La section $s$ vrifie donc les estimes :
\begin{eqnarray*}
\| \nabla^m s \| & \leq & k^{\frac{m}{2}}
\end{eqnarray*}
pour $0 \leq m \leq m_1$ et :
\begin{eqnarray*}
\| \nabla^m \overline{\partial} s \| & \leq & k^{\frac{m}{2}}
\end{eqnarray*}
pour $0 \leq m \leq m_2$.

De mme sur la boule, disons,  $B(y_0,\, 5\times R\times  k^{-\frac{1}{2}})$, les estimes suivantes sont vrifies :
\begin{eqnarray*}
\| \nabla^m c_{y_0}^{-k} \| & \leq & C\, k^{\frac{m}{2}}
\end{eqnarray*}
et :
\begin{eqnarray*}
\| \nabla^m \overline{\partial} c_{y_0}^{-k} \| & \leq & C\, k^{\frac{m}{2}}.
\end{eqnarray*}
(Convenons que, dans cette dmonstration, la valeur du majorant $C$ peut varier d'une ligne  l'autre.)

Posons $f= c_{y_0}^{-k} \times s$. D'aprs la formule de Leibniz, sur la mme boule, cette application $f$ satisfait :
\begin{eqnarray*}
\| f \| & \leq & C
\\ \| \nabla^{m-1} D f \| & \leq & C\, k^{\frac{m}{2}}
\end{eqnarray*}
pour $1 \leq m \leq m_1$ et :
\begin{eqnarray*}
\| \nabla^m \overline{\partial} f \| & \leq &C\, k^{\frac{m}{2}}
\end{eqnarray*}
pour $0 \leq m \leq \mbox{min}\{m_1,\, m_2\}$.

Choisissons un systme de coordonnes complexes centr en le point $y_0$ de $Y$ qui satisfasse les conditions suivantes : 

1. Ce systme est holomorphe en l'origine et isomtrique en l'origine, c'est--dire que sa diffrentielle en $y_0$ est une isomtrie $\C-$linaire entre $T_{y_0}X$ et $\C^n$.

2. Ce systme redresse $Y$ c'est--dire que pour tout point de $X$ proche de $y_0$, il sera quivalent d'appartenir  la sous-varit $Y$ ou de prendre ses coordonnes dans un certain sous-espace vectoriel de $\C^n$, que nous noterons aussi $Y$ par abus de notation.
\\
\\
\indent
(Comme la proposition  dmontrer doit fournir un minorant de $k$, un majorant de $\varepsilon$ et un polynme $P$ qui ne dpendent pas de $y_0$, il faut choisir de telles coordonnes pour chaque point $y_0 \in Y$ de telle sorte que les estimes qui interviendront dans notre dmonstration ne dpendent pas de ce point $y_0$ et que les tailles de leurs domaines de validit ne dpendent pas non plus de $y_0$.  Sans entrer dans les dtails, contentons-nous de dire qu'un tel choix de systmes de coordonnes est possible par compacit de $Y$. ) 
\\
\\
\indent
Le systme de coordonnes centr en $y_0$ permet de dfinir sur les formes diffrentielles et tensorielles la connexion triviale que nous noterons $\nabla_0$.
De mme l'identification locale entre $X$ et $\C^n$ permet de dfinir une structure complexe que nous noterons $J_0$. Enfin, pour tout rel $a> 0$ assez petit, nous noterons $B_0(y_0,\, a)$ la boule de centre $y_0$ et de rayon $a$ pour la distance usuelle de $\C^n$.

Comme $\nabla_0^{m-1}$ est un oprateur diffrentiel d'ordre $m-1$, on a :
\begin{eqnarray*}
\| \nabla_0^{m-1} D f \| & \leq &C \sum_{m'=1}^m \| \nabla^{m'-1} D f \|
\\ & \leq & C\, k^{\frac{m}{2}}
\end{eqnarray*} 
pour $1 \leq m \leq m_1$. De mme :
\begin{eqnarray*}
\| \nabla_0^m \overline{\partial} f \| & \leq &C \sum_{m'=0}^m \| \nabla^{m'} \overline{\partial} f \|
\\ & \leq & C\, k^{\frac{m}{2}}
\end{eqnarray*} 
pour $0 \leq m \leq \mbox{min}\{m_1,\, m_2\}$.

Comme les deux structures presque-complexes concident en $y_0$, elles vrifieront sur la boule $B(y_0,\, 4\times R\times  k^{-\frac{1}{2}})$ :
\begin{eqnarray*}
\| J - J_0\| & \leq & C\, k^{-\frac{1}{2}}
\end{eqnarray*}
Par ailleurs, on a :
\begin{eqnarray*}
\|\nabla_0^m J - \nabla_0^m J_0\| & \leq & C
\end{eqnarray*}
sur la mme boule.

Bilan : en regroupant les deux cas $m=0$ et $m \neq 0$, on obtient :
\begin{eqnarray*}
\|\nabla_0^m J - \nabla_0^m J_0\| & \leq & C\, k^{\frac{1}{2}\mbox{min}\{0,\, m-1\}}
\\ & \leq & C\, k^{\frac{m-1}{2}}.
\end{eqnarray*}
Alors :
\begin{eqnarray*}
\| \nabla_0^m \overline{\partial} f - \nabla_0^m \overline{\partial}_0 f \| & = & \frac{1}{2} \| \nabla_0^m ( Df \circ J) - \nabla_0^m (Df \circ J_0) \| 
\\ & \leq & C \sum_{m'=0}^m \| \nabla_0^{m'} D f \| \times \|\nabla_0^{m-m'} J - \nabla_0^{m-m'} J_0\| 
\\ & \leq & C \sum_{m'=0}^m k^{\frac{m'+1}{2}} \times k^{\frac{m-m'-1}{2}}
\\ & \leq & C\, k^{\frac{m}{2}}
\end{eqnarray*} 
pour $0 \leq m \leq m_1-1$ et donc, par ingalit triangulaire :
\begin{eqnarray*}
\| \nabla_0^m \overline{\partial}_0 f \| & \leq & C\, k^{\frac{m}{2}}
\end{eqnarray*} 
pour $0\leq m \leq  \mbox{min}\{m_1-1,\, m_2\}$.

Comme la diffrentielle du systme de coordonnes est une isomtrie, l'inclusion suivante sera vrifie pour $k$ assez grand :
$$
B_0(y_0,\, 3\times R\times  k^{-\frac{1}{2}}) \subset B(y_0,\, 4\times R\times  k^{-\frac{1}{2}}).
$$
Une consquence vidente est que les estimes que nous avons dmontres sur la seconde boule  seront encore valables sur la premire.

Ces estimes et la proposition (\ref{d bar ech}) impliquent l'existence d'un $g$ holomorphe vrifiant :
\begin{eqnarray*}
\| D^m f -D^m g \| & \leq & C\, k^{\frac{m-1}{2}}
\end{eqnarray*}
pour $0\leq m \leq  \mbox{min}\{m_1-1,\, m_2\}$.
L encore, ces estimes font intervenir {\it a priori} des normes dfinies en utilisant la mtrique riemannienne mais ces normes sont proches de celles qu'on dfinit en utilisant la mtrique usuelle de $\C^n$, pour lesquelles on aura donc des estimes du mme type.

D'aprs le corollaire du thorme (\ref{trans hol}) et la majoration $\|g \| \leq C$, il existe un lment $\alpha \in rL_{y_0}^k$ de norme $\leq \varepsilon$ vrifiant :
\begin{eqnarray*}
\max \left\{ \| g(x) + \alpha \|,\, k^{-\frac{1}{2}} \times \mbox{MS} \left( Dg_Y(x) \right)\right\} & \geq & \frac{\varepsilon}{P\left( \log \frac{1}{\varepsilon} \right)}
\end{eqnarray*}
o $P$ dsigne une fonction polynomiale $>0$. Comme la diffrentielle en $y_0$ du systme de coordonnes est une isomtrie, l'inclusion suivante sera vrifie pour $k$ assez grand :
$$
B(y_0,\, R\times  k^{-\frac{1}{2}}) \subset B_0(y_0,\, 2\times R\times  k^{-\frac{1}{2}}).
$$ 
 L'estime prcdente, valable {\it a priori} sur $Y\cap B_0(y_0,\, 2\times R\times  k^{-\frac{1}{2}})$, sera donc valable sur $Y\cap B(y_0,\, R\times  k^{-\frac{1}{2}})$. 

Les estimes que $g$ vrifie par dfinition impliquent, pour $k$ assez grand :
\begin{eqnarray*}
\max \left\{ \| g(x) - f(x) \|,\, k^{-\frac{1}{2}} \times \left\| Dg(x) - Df(x) \right\| \right\} & \leq & \frac{\varepsilon}{2\, P\left( \log \frac{1}{\varepsilon} \right)}\; \; .
\end{eqnarray*}
(Remarquons que c'est ici que sont intervenues les hypothses $m_1 \geq 2$ et $m_2\ge 1$.)
Par soustraction, on obtient :
\begin{eqnarray*}
\max \left\{ \| f(x) + \alpha \|,\, k^{-\frac{1}{2}} \times \mbox{MS} \left( Df_Y(x)\right)\right\} & \geq & \frac{\varepsilon}{2\, P\left( \log \frac{1}{\varepsilon} \right) }
\end{eqnarray*}
car le module de surjectivit est $1-$lipschitzien. 

Posons :
\begin{eqnarray*}
f' & = & f+\alpha
\\ s' & = & f\, c_{y_0}^k = s+\alpha\, c_{y_0}^k.
\end{eqnarray*}
Alors, en tout point de $B(y_0,\, R\times  k^{-\frac{1}{2}})$ :
\begin{eqnarray*}
\| f' \| & \leq & \| c_{y_0}^{-k} \| \times \| s' \| 
\\ & \leq & C \times  \| s' \|. 
\end{eqnarray*}
La formule de Leibniz :
\begin{eqnarray*}
D f'  & = & c_{y_0}^{-k} \times \nabla s' + s' \times \nabla c_{y_0}^{-k}
\end{eqnarray*}
implique, en utilisant le fait que le module de surjectivit est $1-$lipschitzien, l'estime :
\begin{eqnarray*}
k^{-\frac{1}{2}} \times \mbox{MS} \left( D f_Y' \right) & \leq & 
\| c_{y_0}^{-k}\| \times k^{-\frac{1}{2}} \times \mbox{MS} \left( \nabla s_Y' \right)+ \|s'\| \times k^{-\frac{1}{2}} \times \| \nabla c_{y_0}^{-k}\|
\\ & \leq & C \times \max \left\{ \|s' \| ,\, k^{-\frac{1}{2}} \times \mbox{MS} \left( \nabla s_Y' \right) \right\}
\end{eqnarray*}
valable en tout point de $Y\cap B(y_0,\, R\times  k^{-\frac{1}{2}})$.

En rassemblant les estimes sur $f'$ et sur $\mbox{MS} \left( D f_Y' \right)$ obtenues, on achve la dmonstration de la proposition par le calcul suivant :
\begin{eqnarray*}
\max \left\{ \|s' \| ,\, k^{-\frac{1}{2}} \times \mbox{MS} \left( \nabla s_Y' \right)\right\} 
& \geq & \frac{1}{C} \times \max \left\{ \|f' \| ,\, k^{-\frac{1}{2}} \times \mbox{MS} \left( Df_Y' \right)\right\}
\\ & \geq & \frac{1}{C} \times \frac{\varepsilon}{2\, P\left( \log \frac{1}{\varepsilon} \right) }\; \; .
\end{eqnarray*}

\subsection{Transversalisation sur une partie parpille}\label{trans eparp}
Enonc flottant.
Soient une section $s$ du fibr vectoriel $rL^k$ (avec $k$ assez grand),  $(1,\, 2)-$contrle et $(1,\, 1)-$approximativement holomorphe
et une partie finie $F$ de $Y$ satisfaisant, pour tous $x$, $y\in F$, la condition d'espacement : $d(x,\, y)\geq D \times k^{-\frac{1}{2}}$ 
pour un certain rel $D\geq 1$.

Alors il existe une section $t$ de $rL^k$ de la forme suivante :
 $$
t = \sum_{z \in F} \alpha{^{z}} c_{z}^k
$$
(o, pour tout $z \in F$, le vecteur $\alpha^{z}$ est un lment de $rL_{z}^k$ de norme malor par un $\varepsilon>0$ assez petit),
telle qu'en tout point $y$ de $Y$ vrifiant 
$d(y,\, F)\leq k^{-\frac{1}{2}}$ , le module de transversalit de poids $(1,\, k^{-\frac{1}{2}})$ de $s+t$ le long de $Y$ soit suprieur  $\eta$.
\\
\\
\indent
Attention ! Nous n'affirmons pas que cet nonc soit vrai en toute gnralit. Il ne l'est certainement pas et il nous faut prciser plusieurs choses.
Tout d'abord, comme d'habitude, il nous faut dire que
que \og $\varepsilon$ assez petit \fg\, signifie que $\varepsilon$ est major par un rel $>0$ 
qui dpend des donnes $X$, $L$, $Y$ et $r$
et que \og $k$ assez grand \fg\, signifie que $k$ est suprieur  une borne qui dpend de $\varepsilon$ et
des donnes $X$, $L$, $Y$ et $r$.
Ensuite, plus important, il nous faut dire  quelle condition notre nonc flottant est valable.
La proposition suivante donne une condition suffisante, portant sur les trois rels $\varepsilon>0$, $\eta>0$ et $D\geq 1$,
pour qu'il soit vrifi.
\\
\\
\indent
Proposition.
Il existe une fonction polynomiale $P$  valeurs strictement positives et un rel $C$ (qui dpendent des donnes $X$, $L$, $Y$ et $r$)
telles que l'nonc prcdent sera vrifi si les trois rels $\varepsilon$, $\eta$ et $D$ satisfont les deux conditions :
\begin{eqnarray*}
\frac{\varepsilon}{\eta} & = & P \left( \log \frac{1}{\varepsilon}\right)
\end{eqnarray*}
et :
\begin{eqnarray*}
\frac{\eta}{\varepsilon} & \geq  & C \, \mbox{exp}\left( -\frac{\pi}{4} D^2 \right) .
\end{eqnarray*}
\\
\\
\indent
Dmonstration.
Soient un rel $D \geq 1$ et $F$ une partie finie de $Y$ satisfaisant la condition d'espacement de \og l'nonc flottant \fg.
Pour tout point $z_0 \in F$, la proposition (\ref{trans point}) fournit un lment $\alpha^{z_0} \in rL_{z_0}^k$
de norme $\leq \varepsilon$ tel qu'en tout point de $Y \cap B(z_0, \, k^{-\frac{1}{2}})$, la section 
$s+\alpha^{z_0} c_{z_0}^k$ s'annule transversalement le long de $Y$ avec un module de transversalit de poids $(1,\, k^{-\frac{1}{2}})$ suprieur, 
disons,  $2\eta$ pour un rel $\eta$ vrifiant la premire condition de la proposition. 

Posons :  
$$
t = \sum_{z \in F} \alpha^{z} c_{z}^k.
$$

Soit $y$ un point de $Y$ vrifiant $d(y,\, z_0)\leq k^{-\frac{1}{2}}$ pour un $z_0 \in F$. 
Posons :  
\begin{eqnarray*}
s' & = & s + \alpha^{z_0} c_{z_0}^k
\\ t ' & = & \sum_{z \in F, \, z \neq z_0} \alpha^{z} c_{z}^k
\end{eqnarray*}

L'ingalit triangulaire garantit que tout $z \in F$ diffrent de $z_0$ vrifie $d(y,\, z) \geq (D - 1) \times k^{-\frac{1}{2}}$ et la proposition (\ref{combi}) implique donc les estimes :
\begin{eqnarray*}
\| t'(y) \| & \leq & C \, \varepsilon \, \mbox{exp}\left(-\frac{\pi}{3} (D-1)^2\right)
\\ & \leq & C \, \varepsilon \, \mbox{exp}\left(-\frac{\pi}{4} D^2\right)
\\ \| (\nabla t')_y \| & \leq & C \, \varepsilon \, k^{\frac{1}{2}}\,   \mbox{exp}\left(-\frac{\pi}{3} (D-1)^2\right)
\\ & \leq & C \, \varepsilon \, k^{\frac{1}{2}}\,   \mbox{exp}\left(-\frac{\pi}{4} D^2\right).
\end{eqnarray*}
(o, par convention, la valeur du rel $C$ peut varier d'une ligne  l'autre).

On peut donc fixer une valeur $C>0$ de telle sorte que si $D$ satisfait la seconde condition de la proposition, les deux estimes suivantes seront vrifies :
\begin{eqnarray*}
\| t'(y) \| & \leq & \eta
\\ \| (\nabla t')_y \| & \leq & \eta \times k^{\frac{1}{2}}.
\end{eqnarray*}
Remarquons que les sommes $s+t$ et $s'+t'$ sont gales.
On en dduit, par ingalit triangulaire :
\begin{eqnarray*}
\| (s+t)(y) \| & \geq & \| s'(y) \|-\| t'(y) \|
\\ & \geq & \| s'(y) \|-\eta.
\end{eqnarray*}
Et, le module de surjectivit tant $1-$lischitzien :
\begin{eqnarray*}
\mbox{MS} (\nabla (s_Y +t_Y )_y ) & \geq & \mbox{MS} ((\nabla s'_Y )_y )-\| (\nabla t')_y \|
\\ & \geq & \mbox{MS} ((\nabla s'_Y )_y )-\eta \times k^{\frac{1}{2}}.
\end{eqnarray*}

Les deux estimes prcdentes donnent bien la minoration recherche pour le module de transversalit pondr :
\begin{eqnarray*}
\max \left\{ \| (s+t)(y) \| ,\, k^{-\frac{1}{2}} \times \mbox{MS} (\nabla (s_Y +t_Y )_y ) \right\} & \geq & \max \left\{ \| s'(y) \| ,\, k^{-\frac{1}{2}} \times \mbox{MS} ((\nabla s'_Y )_y ) \right\} - \eta
\\ & \geq & 2\eta - \eta = \eta.
\end{eqnarray*}
(o le terme $2\eta$ provient de la dfinition de $ \alpha {^{z_0}}$). La proposition est dmontre.

\subsection{Nombre de parties espaces ncessaires pour discrtiser $Y$}\label{nombre couleurs}
Proposition.
Soit un rel $D\geq 1$. Alors il existe un entier $n_D$ satisfaisant une majoration du type polynomiale :
$$
n_D \leq P(D)
$$
et tel que, pour tout entier $k\geq 1$, il existe dans $Y$, des parties finies $F_1$, ...., $F_{n_D}$ qui satisfont les deux conditions suivantes :

(i) pour tout $y\in Y$, il existe $z \in \bigcup_{i=1}^{n_D} F_i$ vrifiant : $d(y,\, z)\leq k^{-\frac{1}{2}}$,

(ii) tout entier $i$ compris entre $1$ et $n_D$ et tous points $z$, $z'\in F_i$ vrifient : $d(z,\, z')\geq D\times k^{-\frac{1}{2}}$.  
\\
\\
\indent
Prcisons que $P$ dsigne une fonction polynomiale qui dpend uniquement de $X$ et de $Y$.
\\
\\
\indent
Dmonstration.
Par compacit de $Y$, il existe une partie finie $F$ de $Y$ vrifiant la proprit d'espacement suivante : $d(z,\, z')\geq k^{-\frac{1}{2}}$ et maximale pour cette proprit.
On dfinit alors un graphe abstrait dont les sommets seront les points de $F$ et dont deux sommets $z$ et $z'$ seront relis si et seulement si : $d(z,\, z')< D\times k^{-\frac{1}{2}}$.  
La valence de ce graphe est majore par une fonction polynomiale de $D$ d'aprs la proposition (\ref{taille discr}). La proposition dcoule alors d'un rsultat de thorie des graphes dont la dmonstration est une rcurrence trs lmentaire : le nombre chromatique d'un graphe ne saurait dpasser la valence de plus d'une unit.

\subsection{Un lemme numrique}\label{etoile de mer}
Lemme.
Soient deux rels $0<p<q$ et une suite $(u_n)_{n\geq 1}$ de rels strictement compris entre $0$ et $1$. On suppose vrifie, pour tout $n\geq 2$,
l'ingalit suivante :
\begin{eqnarray}\label{pour u}
u_n &\geq & \frac{u_{n-1}}{\left( \log \frac{1}{u_{n-1}} \right)^p}.
\end{eqnarray}
Posons : 
\begin{eqnarray}
v_n &=&\left( \frac{1}{n}\right)^{nq}.
\end{eqnarray}
Alors il existe un entier $n_0$ tel que tout $n\geq 1$ vrifie :
\begin{eqnarray}\label{lemme concl}
u_n &\geq & v_{n+n_0}.
\end{eqnarray}
\\
\indent
Dmonstration.
L'ingalit $\left( 1-\frac{1}{n}\right)^n \leq e^{-1}$ donne :
\begin{eqnarray}\label{taux}
\frac{v_{n-1}}{v_{n}} & = & (n-1)^q \times \left( 1-\frac{1}{n} \right)^{-nq}
\geq  \left(e\times (n-1)\right)^q.
\end{eqnarray}
Par ailleurs : 
\begin{eqnarray}\label{log trad}
\left( \log \frac{1}{v_{n-1}} \right)^{p} &=&  \left( q \times (n-1)\times \log (n-1)\right)^{p}.
\end{eqnarray}
La comparaison des croissances des membres de droite de (\ref{log trad}) et de (\ref{taux}) montre qu'il existe un entier $n_1$ tel que tout $n\geq n_1$ vrifie :
\begin{eqnarray}
\frac{v_{n-1}}{v_{n}} &\geq & \left(  \log \frac{1}{v_{n-1}} \right)^{p}
\end{eqnarray}
c'est--dire :
\begin{eqnarray}\label{pour v}
v_n &\leq& \frac{v_{n-1}}{\left( \log \frac{1}{v_{n-1}} \right)^p} \;\; .
\end{eqnarray}
Comme la suite $(v_n)$ tend vers $0$, on peut choisir un entier $n_0\geq n_1$ tel que (\ref{lemme concl}) soit vrifie au rang $1$.
Alors pour dmontrer l'ingalit (\ref{lemme concl}) par rcurrence sur $n\geq 1$, il suffit de comparer (\ref{pour u}) et (\ref{pour v}).

\subsection{Conclusion de la dmonstration du thorme de Donaldson-Auroux relatif}

Dans l'nonc du thorme de Donaldson-Auroux relatif, 
on suppose donns un entier $m_{max}$ (disons $\geq 2$) et deux rels $K$ et $\varepsilon>0$. 
On peut supposer $\varepsilon$ arbitrairement petit.
Faisons aussi l'hypothse que $K+\varepsilon$ est major par $1$. 
Cette hypothse n'est pas restrictive car le cas gnral s'y ramne en multipliant les sections donnes par une constante.
Par ailleurs, si le lecteur s'tonne de cette expression $K+\varepsilon$, disons que comme les sections donnes vrifient des estimes faisant intervenir $K$
et qu'on va progressivement les perturber, on verra bien $K+\varepsilon$ intervenir dans les tapes intermdiaires.

La proposition (\ref{combi}) fournit un rel $A>0$ tel que pour tout entier $k\geq 1$, pour toute partie $F$ de $Y$ satisfaisant, pour tous $x$, $y\in F$ la condition
d'espacement : $d(x,y)\geq k^{-\frac{1}{2}}$, si on associe  tout point $y\in F$ un lment $\alpha^y \in rL_y^k$,
et si on pose :
$$
B = \max_{y \in F} \| \alpha^y\|
$$
et :
$$
t=\sum_{y \in F} \alpha^{y} c_yk
$$
les estimes suivantes seront vrifies pour tout entier $m$ compris entre $0$ et $m_{max}$ et pour tout point $x\in X$ : 
$$
\begin{array}{ccl}
\left\| (\nabla^m t)_x \right\| 
& \leq  &
A \times B \times k^{\frac{m}{2}} \times \mbox{exp}(-\frac{k\pi}{3} d^2(x,\, F))
\\
\left\| (\nabla^m \overline{\partial} t)_x \right\| 
& \leq  &
A \times B \times k^{\frac{m}{2}} \times \mbox{exp}(-\frac{k\pi}{3} d^2(x,\, F)).
\end{array}
$$
On peut supposer : $\frac{\varepsilon}{A} <\frac{1}{2e}$, quitte  augmenter $A$.
Notons $C$ et $P$ un rel et un polynme fournis par la proposition (\ref{trans eparp}).
Posons $\varepsilon_1=\frac{\varepsilon}{A}$ et
$\eta_1 =\frac{\varepsilon_1}{P \left( \log \frac{1}{\varepsilon_1} \right)}$.
Puis dfinissons, pour tout entier $i\geq 2$, les rels strictement positifs $\varepsilon_i$ et $\eta_i$ par les relations de rcurrence :
$$\varepsilon_{i+1} = \min \left( \varepsilon_1\, ,\;  \frac{\eta_i}{2A}\right) \;\; \mbox{ et }\;\; \eta_{i+1} = \min \left(\frac{\eta_i}{2}\, ,\; \frac{\varepsilon_{i+1}}{P \left( \log \frac{1}{\varepsilon_{i+1}} \right)} \right).$$

Nous aurons besoin du rsultat suivant :
\\
\\
\indent
Lemme.
Il existe un rel $D\geq 1$ tel que l'ingalit suivante soit satisfaite : 
$$\frac{\eta_i}{\varepsilon_i} \geq C\, \mbox{exp}(-\frac{\pi}{4} D^2)$$
pour tout entier $i$ compris entre $1$ et l'entier $n_D$ fourni par la proposition (\ref{nombre couleurs}).
\\
\\
\indent
Dmonstration du lemme.
Vrifions qu'il existe un rel $p$ vrifiant pour tout $i$ l'ingalit :
$$
\varepsilon_{i+1}\geq \frac{\varepsilon_i}{\left( \log \frac{1}{\varepsilon_1} \right)^p}.
$$
Par dfinition de $\varepsilon_{i+1}$, il suffit de le vrifier dans les deux cas suivants :

\noindent {\it 1er cas : }$\varepsilon_{i+1} = \varepsilon_1$. 

Le rsultat est vident.

\noindent {\it 2nd cas : }$\varepsilon_{i+1} = \frac{\eta_i}{2A}$. 

La dfinition de $\eta_i$ nous conduit  distinguer deux sous-cas :

{\it 1er sous-cas : }$\eta_{i} =  \frac{\eta_{i-1}}{2}$. 

Alors on a :
$
\varepsilon_{i+1} = \frac{\eta_{i-1}}{4A}\geq \frac{\varepsilon_i}{2}
$
par dfinition de $\varepsilon_i$. Le rsultat en dcoule.

{\it 2nd sous-cas : }$\eta_{i} =  \frac{\varepsilon_{i}}{P \left( \log \frac{1}{\varepsilon_{i}} \right)}$. 

Alors on a :
$
\varepsilon_{i+1} =\frac{1}{2A}\times \frac{\varepsilon_{i}}{P \left( \log \frac{1}{\varepsilon_{i}} \right)}
$
et le rsultat en dcoule.
\\
\\
\indent
Ces vrifications tant faites, le lemme (\ref{etoile de mer}) implique l'existence d'un entier $i_0$ tel que tout $i\geq 1$ satisfasse :
$$
\varepsilon_i \geq \left( \frac{1}{i+i_0} \right) ^{q\times  (i+i_0)}
$$
pour, disons, $q=p+1$.
Les entiers $i\leq n_D$ vrifient donc :
$$
\varepsilon_i \geq \left( \frac{1}{n_D+i_0} \right) ^{q\times (n_D+i_0)}.
$$
Le membre de droite est minor par l'exponentielle d'une fonction polynomiale de $n_D$
qui est lui-mme polynomial en $D$. 
Il existe donc une fonction polynomiale $g$ vrifiant :
$$
\log \frac{1}{\varepsilon_i} \leq g(D).
$$
Par dfinition de $\eta_i$, on en dduit l'existence d'une fonction polynomiale $h$  valeurs $>0$ vrifiant :
$$
\frac{\eta_i}{\varepsilon_i} \geq \frac{1}{h(D)}.
$$
Pour $D$ assez grand, ce rapport est bien minor par $C\, \mbox{exp}(-\frac{\pi}{4} D^2)$, ce qui achve la dmonstration du lemme. 
\\
\\
\indent
Poursuivons la dmonstration du thorme de Donaldson-Auroux relatif.
Dans la suite, nous noterons  $D$ le rel obtenu en appliquant ce lemme
et nous noterons bien sr $n_D$ l'entier correspondant fourni par la proposition (\ref{nombre couleurs}).
Soient des parties $F_1$, ...., $F_{n_D}$ de $Y$, qui satisfont les conditions (i) et (ii) de la proposition (\ref{nombre couleurs}).
Pour tout $i$ compris entre $0$ et $n_D$, notons $G_i$ la runion de $F_1$, $F_2$, ..., $F_i$.

Un second lemme est ncessaire, qui transversalise une section sur certaines parties de $Y$ :
\\
\\
\indent
Lemme.
Soient un entier $i$ compris entre $1$ et $n_D$ et une section $s$ de $rL^k$ (avec $k$ assez grand), qu'on suppose $(K, 2)-$contrle et $(K,1)-$approximativement holomorphe.
Pour toute partie $P$ de $X$,  nous noterons $V_k(P)$ l'ensemble des points de $X$ dont la distance  $P$ est majore par $k^{-\frac{1}{2}}$.
Il existe une section $t_i$ de $rL^k$ satisfaisant les deux conditions suivantes :

(i) En tout point de $Y \cap V_k (G_i)$, le module de transversalit de poids $(1,\, k^{-\frac{1}{2}})$ de $s+t_i$ le long de $Y$ est suprieur  $\eta_i$.

(ii) La section $t$ est de la forme :
$$
t_i=\sum_{y \in G_i}  \alpha^{y} c_y^k
$$
o, pour tout $j$ compris entre $1$ et $i$ et tout $y\in F_j$, on dsigne par $\alpha^{y}$ un lment de $rL_y^k$ de norme $\leq \varepsilon_j$. 

Prcisons que \og $k$ assez grand \fg\, signifie que $k$ est suprieur  une borne qui dpend 
des donnes $X$, $L$, $Y$, $r$, $K$ et $\varepsilon$ ainsi que des choix de $C$, $P$, $D$ et $n_D$.
\\
\\
\indent
Dmonstration du lemme, par rcurrence sur $i$.

Soit un entier $i$ compris entre $0$ et $n_D -1$. 
On fait
l'hypothse de rcurrence qu'il existe une section $t_i$ vrifiant les conditions (i) et (ii)
(sauf si $i$ est nul, auquel cas on pose simplement $t_i=0$).

Remarquons que la somme $s+t_i$ est $(1, 2)-$contrle et $(1,1)-$approximativement holomorphe. 
Le prcdent lemme et la proposition (\ref{trans eparp}) impliquent donc l'existence d'une section $u_{i+1}$ de $rL^k$ satisfaisant les deux conditions suivantes :

(a) La section $u_{i+1}$ est de la forme :
$$
u_{i+1}=\sum_{y \in F_{i+1}} \alpha^{y} c_y^k
$$
o, pour tout $y\in F_{i+1}$, on dsigne par $\alpha^{y}$ un lment de $rL_y^k$ de norme $\leq \varepsilon_i$.

(b) En tout point de $Y \cap V_k (F_{i+1})$, le module de transversalit de poids $(1,\, k^{-\frac{1}{2}})$ de $s+t_i+u_{i+1}$ le long de $Y$ est suprieur  
$\frac{\varepsilon_{i+1}}{P \left( \log \frac{1}{\varepsilon_{i+1}} \right)} $ et donc  $\eta_{i+1}$.
\\
\\
\indent
Posons : $t_{i+1} = u_{i+1} + t_i$.
Comme le module de surjectivit est $1-$lipschitzien, la section $s+t_{i+1}$ vrifie, en tout point de $Y \cap V_k(G_i)$, la minoration suivante :
\begin{eqnarray*}
\max \left\{\| s+t_{i+1} \| ,\, k^{-\frac{1}{2}} \times \mbox{MS} \left( \nabla (s+t_{i+1} )_Y\right) \right\} 
& \geq & \max \left\{\| s+t_i \| ,\,  k^{-\frac{1}{2}} \times \mbox{MS} \left( \nabla (s+t_i )_Y \right) \right\}  
\\ & & - \; \; \max \left\{ \| u_{i+1} \| ,\,  k^{-\frac{1}{2}} \times \| \nabla u_{i+1} \| \right\} 
\\ & \geq & \eta_i - A \times  \varepsilon_{i+1} 
\\ & \geq &\frac{\eta_i}{2} \; \geq \;  \eta_{i+1}.
\end{eqnarray*}
Comme $V_k (G_{i+1})$ est la runion de $V_k(G_i)$ et de $V_k (F_{i+1})$, la section $t_{i+1}$ satisfait bien les conditions (i) et (ii) au rang $i+1$.
Par rcurrence, le lemme est dmontr.
\\
\\
\indent
Suite et fin de la dmonstration du thorme.
Le lemme (pour $i=n_D$) fournit une section $t$ telle que la somme $s+t$ ait un module de transversalit pondr le long de $Y$ suprieur  $\eta_{n_D}$.
Par ailleurs, cette section $t$ est bien $(\varepsilon,m_{max} )-$contrle et $(\varepsilon,m_{max} )-$approximativement holomorphe d'aprs
l'ingalit : $$\max_{1\leq i \leq n_D} \varepsilon_i \leq \frac{\varepsilon}{A}$$ et la dfinition de $A$.

\vskip13pt
\footnotesize{
\indent 

CMI, LATP, UNIVERSITE D'AIX-MARSEILLE, 
39 RUE F. JOLIOT-CURIE, 
F-13453 MARSEILLE CEDEX 13 FRANCE
\\
{\it \noindent E-mail address :} jean-paul.mohsen@univ-amu.fr}

\end{document}